\def\TheMagstep{\magstep1}	
\def\PaperSize{letter}		

\magnification=\magstep1

\let\:=\colon  
\let\vf=\varphi  \let\o=\circ \let\?=\overline

\let\Sum=\sum \def\sum{\Sum\nolimits}

\def\IC{{\bf C}} 
\def\IP{{\bf P}} \def\mY{{\bf m}_Y}

\def\Pd #1#2{{\partial#1\over\partial#2}}
\def\sdep #1{#1^\dagger}	
\def\and{\hbox{ and }}
\def\Wf{\hbox{\rm W$_f$\hskip 2pt}}
\def\WF{\hbox{\rm W$_F$\hskip 2pt}}	
\def\Af{\hbox{\rm A$_f$\hskip 2pt}}
\def\Al{\hbox{\rm A$_l$\hskip 2pt}}
\def\AF{\hbox{\rm A$_F$\hskip 2pt}}

\def\DONE{*!*}
\def\NextDef #1 {\def\NextOne{#1}%
 \ifx\NextOne\DONE\let\next\relax
 \else\expandafter\xdef\csname#1\endcsname{\TheOp}
  \let\next\NextDef
 \fi \next}
\def\TheOp{\mathop{\rm\NextOne}}
 \NextDef 
  Projan Supp Proj Sym Spec Hom cod Ker dist
 *!*
\def\TheOp{{\cal\NextOne}}
\NextDef 
  E F G H I J M N O R S
 *!*
\def\TheOp{\hbox{\rm\NextOne}}
\NextDef 
 A ICIS 
 *!*

\def\item#1 {\par\indent\indent\indent\indent \hangindent4\parindent
 \llap{\rm (#1)\enspace}\ignorespaces}
 \def\inpart#1 {{\rm (#1)\enspace}\ignorespaces}
 \def\part {\par\inpart}

\catcode`\@=11		

\def\vfootnote#1{\insert\footins\bgroup
 \eightpoint 
 \interlinepenalty\interfootnotelinepenalty
  \splittopskip\ht\strutbox 
  \splitmaxdepth\dp\strutbox \floatingpenalty\@MM
  \leftskip\z@skip \rightskip\z@skip \spaceskip\z@skip \xspaceskip\z@skip
  \textindent{#1}\footstrut\futurelet\next\fo@t}

\def\p.{p.\penalty\@M \thinspace}
\def\pp.{pp.\penalty\@M \thinspace}
\def\(#1){{\rm(#1)}}\let\leftp=(
\def\activeleftp{\catcode`\(=\active}
{\activeleftp\gdef({\ifmmode\let\next=\leftp \else\let\next=\(\fi\next}}

\def\sct#1\par
  {\removelastskip\vskip0pt plus2\normalbaselineskip \penalty-250 
  \vskip0pt plus-2\normalbaselineskip \bigskip
  \centerline{\smc #1}\medskip}

\newcount\sctno \sctno=0
\def\sctn{\advance\sctno by 1 
 \sct\number\sctno.\quad\ignorespaces}

\def\dno#1${\eqno\hbox{\rm(\number\sctno.#1)}$}
\def\Cs#1){\unskip~{\rm(\number\sctno.#1)}}

\def\proclaim#1 #2 {\medbreak
  {\bf#1 (\number\sctno.#2)}\enspace\bgroup\activeleftp
\it}
\def\endproclaim{\par\egroup\medskip}
\def\pf{\endproclaim{\bf Proof.}\enspace}
\def\lem{\proclaim Lemma } \def\prp{\proclaim Proposition }
\def\cor{\proclaim Corollary }	\def\thm{\proclaim Theorem }
\def\rmk#1 {\medbreak {\bf Remark (\number\sctno.#1)}\enspace}
\def\eg#1 {\medbreak {\bf Example (\number\sctno.#1)}\enspace}

\parskip=0pt plus 1.75pt \parindent10pt
\hsize29pc
\vsize44pc
\abovedisplayskip6pt plus6pt minus2pt
\belowdisplayskip6pt plus6pt minus3pt

\def\TRUE{TRUE}	
\ifx\DoublepageOutput\TRUE \def\TheMagstep{\magstep0} \fi
\mag=\TheMagstep

\newskip\vadjustskip \vadjustskip=0.5\normalbaselineskip
\def\centertext
 {\hoffset=\pgwidth \advance\hoffset-\hsize
  \advance\hoffset-2truein \divide\hoffset by 2\relax
  \voffset=\pgheight \advance\voffset-\vsize
  \advance\voffset-2truein \divide\voffset by 2\relax
  \advance\voffset\vadjustskip
 }
\newdimen\pgwidth\newdimen\pgheight
\def\letter{letter}\def\AFour{AFour}
\ifx\PaperSize\letter
 \pgwidth=8.5truein \pgheight=11truein 
 \message{- Got a paper size of letter.  }\centertext 
\fi
\ifx\PaperSize\AFour
 \pgwidth=210truemm \pgheight=297truemm 
 \message{- Got a paper size of AFour.  }\centertext
\fi

 \newdimen\fullhsize \newbox\leftcolumn
 \def\fulline{\hbox to \fullhsize}
\def\doublepageoutput
{\let\lr=L
 \output={\if L\lr
          \global\setbox\leftcolumn=\columnbox \global\let\lr=R%
        \else \doubleformat \global\let\lr=L\fi
        \ifnum\outputpenalty>-20000 \else\dosupereject\fi}%
 \def\doubleformat{\shipout\vbox{%
        \fulline{\hfil\hfil\box\leftcolumn\hfil\columnbox\hfil\hfil}%
				}%
		  }%
 \def\columnbox{\vbox
   {\makeheadline\pagebody\makefootline\advancepageno}%
   }%
 \fullhsize=\pgheight \hoffset=-1truein
 \voffset=\pgwidth \advance\voffset-\vsize
  \advance\voffset-2truein \divide\voffset by 2
  \advance\voffset\vadjustskip
 \let\firstheadline=\hfil
 
}
\ifx\DoublepageOutput\TRUE \doublepageoutput \fi

 \font\twelvebf=cmbx12		
 \font\smc=cmcsc10		

\def\eightpoint{\eightpointfonts
 \setbox\strutbox\hbox{\vrule height7\p@ depth2\p@ width\z@}%
 \eightpointparameters\eightpointfamilies
 \normalbaselines\rm
 }
\def\eightpointparameters{%
 \normalbaselineskip9\p@
 \abovedisplayskip9\p@ plus2.4\p@ minus6.2\p@
 \belowdisplayskip9\p@ plus2.4\p@ minus6.2\p@
 \abovedisplayshortskip\z@ plus2.4\p@
 \belowdisplayshortskip5.6\p@ plus2.4\p@ minus3.2\p@
 }
\newfam\smcfam
\def\eightpointfonts{%
 \font\eightrm=cmr8 \font\sixrm=cmr6
 \font\eightbf=cmbx8 \font\sixbf=cmbx6
 \font\eightit=cmti8 
 \font\eightsmc=cmcsc8
 \font\eighti=cmmi8 \font\sixi=cmmi6
 \font\eightsy=cmsy8 \font\sixsy=cmsy6
 \font\eightsl=cmsl8 \font\eighttt=cmtt8}
\def\eightpointfamilies{%
 \textfont\z@\eightrm \scriptfont\z@\sixrm  \scriptscriptfont\z@\fiverm
 \textfont\@ne\eighti \scriptfont\@ne\sixi  \scriptscriptfont\@ne\fivei
 \textfont\tw@\eightsy \scriptfont\tw@\sixsy \scriptscriptfont\tw@\fivesy
 \textfont\thr@@\tenex \scriptfont\thr@@\tenex\scriptscriptfont\thr@@\tenex
 \textfont\itfam\eightit	\def\it{\fam\itfam\eightit}%
 \textfont\slfam\eightsl	\def\sl{\fam\slfam\eightsl}%
 \textfont\ttfam\eighttt	\def\tt{\fam\ttfam\eighttt}%
 \textfont\smcfam\eightsmc	\def\smc{\fam\smcfam\eightsmc}%
 \textfont\bffam\eightbf \scriptfont\bffam\sixbf
   \scriptscriptfont\bffam\fivebf	\def\bf{\fam\bffam\eightbf}%
 \def\rm{\fam0\eightrm}%
 }

\def\today{\ifcase\month\or	
 January\or February\or March\or April\or May\or June\or
 July\or August\or September\or October\or November\or December\fi
 \space\number\day, \number\year}
\nopagenumbers
\headline={%
  \ifnum\pageno=1\firstheadline
  \else
    \ifodd\pageno\oddheadline
    \else\evenheadline\fi
  \fi
}
\let\firstheadline\hfill
\def\oddheadline{
 \hfil\headtitle\hfil\llap{\folio}}
\def\evenheadline{\eightpoint\rlap{\folio}
 \hfil\author\hfil
}
\def\headtitle{\title}

 \newcount\refno \refno=0	 \def\NoKey{*!*}
 \def\MakeKey{\advance\refno by 1 \expandafter\xdef
  \csname\TheKey\endcsname{{\number\refno}}\NextKey}
 \def\NextKey#1 {\def\TheKey{#1}\ifx\TheKey\NoKey\let\next\relax
  \else\let\next\MakeKey \fi \next}
 \def\RefKeys #1\endRefKeys{\expandafter\NextKey #1 *!* }
\def\SetRef#1 #2,#3\par{%
 \hang\llap{[\csname#1\endcsname]\enspace}%
  \ignorespaces{\smc #2,}
  \ignorespaces#3\unskip.\endgraf
 }
 \newbox\keybox \setbox\keybox=\hbox{[8]\enspace}
 \newdimen\keyindent \keyindent=\wd\keybox
\def\references{
  \bgroup   \frenchspacing   \eightpoint
   \parindent=\keyindent  \parskip=\smallskipamount
   \everypar={\SetRef}}
\def\endreferences{\egroup}

 \def\serial#1#2{\expandafter\def\csname#1\endcsname ##1 ##2 ##3
  {\unskip\ #2 {\bf##1} (##2), ##3}}
 \serial{ajm}{Amer. J. Math.}
  \serial {aif} {Ann. Inst. Fourier}
 \serial{asens}{Ann. Scient. \'Ec. Norm. Sup.}
 \serial{comp}{Compositio Math.}
 \serial{conm}{Contemp. Math.}
 \serial{crasp}{C. R. Acad. Sci. Paris}
 \serial{dlnpam}{Dekker Lecture Notes in Pure and Applied Math.}
 \serial{faa}{Funct. Anal. Appl.}
 \serial{invent}{Invent. Math.}
 \serial{ma}{Math. Ann.}
 \serial{mpcps}{Math. Proc. Camb. Phil. Soc.}
 \serial{ja}{J. Algebra}
 \serial{splm}{Springer Lecture Notes in Math.}
 \serial{tams}{Trans. Amer. Math. Soc.}

\def\UThin{\penalty\@M \thinspace\ignorespaces}
\def\relaxnext@{\let\next\relax}
\def\cite#1{\relaxnext@
 \def\nextiii@##1,##2\end@{\unskip\space{\rm[\SetKey{##1},\let~=\UThin##2]}}%
 \in@,{#1}\ifin@\def\next{\nextiii@#1\end@}\else
 \def\next{{\rm[\SetKey{#1}]}}\fi\next}
\newif\ifin@
\def\in@#1#2{\def\in@@##1#1##2##3\in@@
 {\ifx\in@##2\in@false\else\in@true\fi}%
 \in@@#2#1\in@\in@@}
\def\SetKey#1{{\bf\csname#1\endcsname}}

\catcode`\@=12 
\def\title{ The multiplicity of pairs of modules and hypersurface 
singularities}
\def\author{Terence Gaffney}
\RefKeys BMM BS B-R Bo1 F G-1 G-2 G-3 G-4 G-5 G-6 G-7 G-G GK GK2 GM Gas H-M J 
K-T KT1 Le-T LJT Ma M N P P1 R S  T-1 T-2 Z Z1
 \endRefKeys

\def\topstuff{\leavevmode
 \bigskip\bigskip
 \centerline{\twelvebf \title}
 \bigskip
 \centerline{\author}
\bigskip\bigskip}
\topstuff
\sct Introduction

In \cite{G-5} the author introduced the notion of the multiplicity of a 
pair of modules as a way of working with modules of non-finite colength.
Some applications of this notion to equisingularity problems were 
described in \cite{G-6}.  The invariants introduced using
  this tool have the advantage that they must be independent of the
parameters in the family when the stratification condition they describe
holds.  These invariants provide a framework for studying the
equisingularity conditions $W$, $W_f$ and $A_f$ for very general
families
of spaces and functions. In this paper we will illustrate the use of
these
invariants in the study of families of functions with non-isolated
singularities and show how the invariants arise naturally in the work of
Pellikaan (\cite{P}, \cite{P1})and Zaharia (\cite{Z}, \cite{Z1}). 
Pellikaan studied functions $f$ whose singular set 
was an isolated complete intersection singularity (ICIS) of dimension 1, 
Zaharia those of of dimension 2. 

The principal tool for connecting the multiplicity of the pair with 
geometry is the multiplicity polar theorem (Theorem 1.1) which we review 
in section 1. This theorem is used to relate multiplicity information at the special fiber of a family with information at the generic fiber.
As an illustration of the theorem we use it to give a geometric 
interpretation of the multiplicity of a module (Theorem 1.2). This 
interpretation is then used in
remark 1.3 to connect the multiplicity of a module with Fulton's k-th 
degeneracy class. In section 2 we show how the multiplicity
of the  pair $(J(f),I)$ appears naturally in the work of Pellikaan and 
give two formulas for it. The first formula relates this multiplicity
to the number of D$_{\infty}$ and A$_1$ points appearing in a deformation 
of $f$. The second formula shows that the multiplicity of the pair 
   $(J(f),I)$ if defined is actually the length of $I/J(f)$. This length 
is Pellikaan's invariant $j(f)$. Both formulas are contained in Theorem 
2.3 and its proof. These formulas are
used to give a new formula for  the L\^e number of dimension 0 (Proposition 2.4). 
( Cf. \cite{Ma} for details on the L\^e numbers.)

In section 3, we extend the results of Pellikaan for singular sets of 
dimension 1 to ICIS of dimension d, then use these results to prove 
extensions of the theorems of section 2. The computation of the formula 
for 
 the L\^e number of dimension 0 uses Zaharia's computation of the homology 
of the Milnor fiber.

These formulae suggest in general that the L\^e number of dimension 0 is 
the sum of the invariant which controls the \Af condition, (which in turn is related the multiplicity of a pair of
ideals), and invariants of dimension 0 related to the other singularity 
types in the singular set of $f$.

Section 3 also shows that the condition that $j(f)$ is finite imposes 
strong 
restrictions on $f$--there must exist a set of generators of $I$, 
$\{g_1,\dots,g_p\}$ such that
$f=\sum g^2_i$. This implies that every such function is the composition 
of a function $h$ with a Morse singularity at the origin and the map $G$ 
whose components are generators
of $I$. In particular, all of the germs of type $D(d,p)$, with $d>1$, studied by 
Pellikaan have $j(f)=\infty$, contrary to assertions made in Remark 5.3 on 
page 52 of \cite{P}
and in Remark 5.4 on page 373 of \cite{P1}.

In section 4  we then use the multiplicity
of the  pair to give a necessary and sufficent condition for the \Af 
condition to hold for a family of functions $f_y$ (Theorem 4.5). 
The proof of this result involves a new trick which is used to pass 
information from strata in the singular set of $f$ to the ambient geometry 
of 
$f_0$. This enables us to drop the hypothesis that the ``natural" 
stratification of the singular set of $f$ satisfies Whitney A.

In the case that the singular locus of $f_0$ is an ICIS of dimension 1, we 
use the relation between our invariant and the L\^e numbers,
 to show that a strong form of the \AF condition also 
implies that the L\^e numbers are constant as well (Corollary 4.7). This 
is used to show that in this situation the strong form of the \Af 
condition implies the triviality of the Milnor
fibrations (Corollary 4.8). In example 4.9, by modifying the example of 
Briancon-Speder we show that 
both the \Af condition and topological triviality of the family may hold, 
yet the L\^e numbers may not be constant. It remains open whether the 
strong form of the \Af condition implies
the L\^e number of dimension $0$ is constant in general, or if the strong form of \Af is 
needed if the dimension of $S(f)=1$.

We then discuss the \Wf condition for the situation of Theorem 4.5. Here 
we show that the independence from parameter of a single invariant is all that is required 
for a \Wf-Whitney 
stratification of a family of functions, which implies the topological 
triviality of the family (Theorem 4.10). 
This invariant is then related to the relative polar
multiplicites of the  members of the family and the multiplicity of the 
pair that is used to control the \Af condition (Corollary 4.13).
 In turn, this implies that the 
\Af condition combined with the independence from parameter of the 
relative polar multiplicities implies
 that we have a \Wf-Whitney stratification (Corollary 4.14).

The application of the multiplicity of the pair to equisingularity 
problems grew out of a long series of conversations with Steven Kleiman; 
the author thanks him for his
encouragement. The author also thanks David Massey and James Damon for 
helpful conversations, and the referee for his careful reading of the paper, and helpful suggestions.

\sctn The multiplicity polar theorem

In this paper we work with complex analytic sets and maps. Let $\O_X$ 
denote the strucure sheaf on a complex analytic space $X$.
If a module $M$ has finite colength in $\O^p_{X,x}$, it is possible to 
attach a number 
to the module, its Buchsbaum-Rim multiplicity (\cite{B-R}). We can also 
define the multiplicity
of a pair of modules $M\subset N$, $M$ of finite colength in $N$, as well, 
even if $N$ does not have finite colength in $\O^p_X$. 
We recall how to construct these numbers following the approach
 of Kleiman and Thorup (\cite{K-T}).

Given a submodule $M$ of a free $\O_X$ module $F$ of rank $p$, we can 
associate a subalgebra ${\cal R }(M)$ of the symmetric
$\O_X$ algebra on $p$ generators. This is known as the Rees algebra of 
$M$. If $(m_1,\dots,m_p)$ is an element of $M$
 then $\Sum m_iT_i$ is the corresponding element of 
${\cal R}(M)$. Then $\Projan ({\cal R} (M))$, the projective analytic 
spectrum of ${\cal R}(M)$ is the closure of 
the projectivised row spaces of $M$ at points where the rank
of a matrix of generators of $M$ is maximal. Denote the projection to $X$ 
by $c$, or by $c_M$ where there is ambiguity.

 If $M$ is a submodule of $N$ or $h$ is a section of $N$, then
 $h$ and  $M$  generate ideals on $\Projan{\cal R}(N)$; denote them
by $\rho(h)$ and $\rho({\cal M})$. If we can express $h$ in terms of a set 
of generators $\{n_i\}$ of $N$ as $\Sum g_in_i$, 
then in the chart in which $T_1\ne 0$, we can express a generator of 
$\rho(h)$ by $\Sum g_iT_i/T_1$. 
Having defined the ideal sheaf
$\rho({\cal M})$, we blow up by it.

On the blowup $B_{\rho({\cal M})}(\Projan {\cal R}(N))$ we have two 
tautological bundles, one the pullback of 
the bundle on $\Projan{\cal R}(N)$,
 the other coming from $\Projan{\cal R}(M)$; denote the corresponding 
Chern classes by $l_M$ and $l_N$, and denote the exceptional divisor by
$ D_{M,N}$. 
Suppose the generic 
rank of $N$ (and hence of $M$) is $e$.
 Then 
the multiplicity of a pair of modules $M,N$ is:

$$e(M,N)=\Sum_{j=0}^{d+e-2} \int D_{M,N} \cdot l^{d+e-2-j}_{M}\cdot 
l^j_{N}.$$

 The multiplicity of the pair is well defined as long as the set of points where $N$ is not integrally dependent on $M$ is isolated (\cite{K-T}). If the pair is $M$ and  $\O^p_X$, then this condition implies
that $M= \O^p_X$ except at isolated points so $M\subset \O^p_X$ is of finite colength and  the multiplicity of $M$ 
is the multiplicity of the pair $(M,\O^p_X)$.
 Later in this section we 
will give a new geometric interpretation of this number based on polar 
methods.

If $\O_{X^d,x}$ is Cohen-Macauley, and $M$ has $d+p-1$ generators then 
there is a useful relation between $M$ and its ideal of maximal minors; 
the multiplicity of $M$ is the colength of $M$, is the
 colength of the ideal of maximal minors,
 by some theorems of Buchsbaum and Rim \cite{B-R}, 2.4 \p.207, 4.3 and
4.5 \p.223. In section 2 we will see a first generalization of this result 
to pairs of modules.

We next develop the notion of polar varieties which is the other term in 
the multiplicity polar theorem.

Assume we have a module $M$ which is a submodule of a free module on $X^d$, 
an equidimensional,  analytic space, reduced off a nowhere dense subset of 
$X$,   and the generic rank of $M$ is $e$ on each 
component of $X$. The hypothesis on the equidimensionality of $X$ and on 
the rank of $M$ ensures that $\Projan{\cal R}(M)$ 
is equidimensional of 
dimension $d+e-1$. Note that $\Projan{\cal R}(M)$ can be embedded in 
$X\times \IP^{r-1}$, 
provided we can chose a set of generators of $M$ 
with $r$ elements.
The {\it polar variety of codimension $k$} of $M$ in $X$ denoted 
$\Gamma_k(M)$ is constructed by intersecting $\Projan{\cal R}(M)$
 with $X\times H_{e+k-1}$ where 
$H_{e+k-1}$ is a general plane of codimension $e+k-1$, then projecting to 
$X$. This notion was developed by Teissier 
in the case where 
$M=JM(F)$, $X=F^{-1}(0)$ (\cite{T-2}). Think of $H$ as the projectivised 
row space of a linear submersion $\pi$. 
 Then $\Gamma_k(JM(F))$ consists of
 the set 
of points where the matrix formed from $D\pi$ and $DF$ has less than 
maximal rank, hence greater than minimal kernel
 rank. These
are the points where
the restriction to $X$ of $\pi$ is singular. In general, think of 
$\Gamma_k(M)$ as the set of points 
where the module whose matrix of generators consists of the matrix of 
generators of $M$ augmented by the rows of 
the linear submersion
$\pi$, has less than maximal rank $n-k+1$. 
When we consider $M$ as part of a pair of modules $M,N$, where the generic 
rank of $M$ is the same as the generic rank of $N$,
 then other polar varieties become interesting as well. In brief, we can 
intersect
 $B_{\rho({\cal M})}(\Projan {\cal R}(N))\subset X\times 
\IP^{N-1}\times\IP^{p-1}$ with a mixture of hyperplanes from the two 
projective spaces
which are factors of the space in which the blowup is embedded. We can 
then push these intersections down to 
$\Projan {\cal R}(N)$ or
$X$ as is convenient, getting mixed polar varieties in $\Projan {\cal 
R}(N)$ or in $X$. These mixed varieties play an 
important role in the proof of the multiplicity-polar theorem, the theorem 
we next describe.

Setup: We suppose we have families  of modules $M\subset  N$, $M$ and $N$
 submodules of a free module $F$ of rank $p$
 on an equidimensional family of spaces with equidimensional
 fibers ${\cal X}^{d+k}$, ${\cal X}$ a family over a smooth base
$Y^k$. We assume that the generic rank of $M$, $N$ is $e\le p$.  Let 
$P(M)$ denote $\Projan {\cal R}(M)$, $\pi_M$
 the projection to ${\cal X}$.  
let $C(M)$ denote the locus of points where $M$ is not free, ie. the 
points where the rank of $M$ is less 
than $e$, $C(\Projan {\cal R}(M))$
its inverse image under $\pi_M$, $C({\cal M})$ the cosupport of 
$\rho({\cal M})$ in $P(\Projan {\cal R}(N))$.

We will be interested in computing the change in the multiplicity of the 
pair $(M,N)$, denoted  $\Delta (e(M,N))$. We will
assume that the integral closures of $M$ and $N$ agree  off a set $C$ of 
dimension $k$ which is finite over $Y$, and assume
we are working on a sufficently small neighborhood of the origin, that 
every component of $C$ contains the origin in its
closure. Then $e(M,N,y)$ is the sum of the multiplicities of the pair at 
all points in the fiber of $C$ over $y$, and $\Delta (e(M,N))$
is the change in this number from $0$ to a generic value of $y$. If we 
have a set $S$ which is finite over $Y$, then we can project $S$ to $Y$, 
and the degree
of the branched cover at $0$ is ${\rm  mult}_y S$. (Of course, this is 
just the number of points in the fiber of $S$ over our generic $y$.)

We can now state our theorem.

\thm 1 Suppose in the above setup we have that $\?M=\?N$ off a set $C$ of 
dimension $k$ which is finite over $Y$. 
Suppose further
that $C(\Projan {\cal R}(M))(0)=C(\Projan {\cal R}(M(0)))$ except possiby 
at the points which project to $0\in {\cal X}(0)$. 
Then, for $y$ a 
generic point of $Y$,

$$\Delta (e(M,N))={\rm  mult}_y \Gamma_{d}(M)-{\rm mult}_y 
\Gamma_{d}(N).$$

\pf The proof in the ideal case appears in \cite{G-6}; the general proof 
will appear in \cite{G-7}.

Now we describe an application of the result to the simple case where $N$ 
is free.

The following geometric interpretation of the multiplicity of an ideal is 
well known. 
Given an ideal $I$ of finite colength in $ \O_{X,x}$, $X^d$ 
equidimensional, choose $d$ elements $(f_1,\dots,f_d)$ of 
$I$ which generate a reduction of $I$. (Recall that if $M$ is a submodule 
of $N$, 
then $M$ is {\it reduction} of $N$  if they have the same integral 
closure.) Then the multiplicity of $I$ is the degree at $x$ of $F$ where $F$
is the branched cover defined by the map-germ with components  
$(f_1,\dots,f_d)$, 
or the number of points in a fiber of $F$ over a regular value close to 
$0$.

We wish to give a similar interpretation of the multiplicity of a module. 

\thm 2 Given $M$ a submodule of  $\O^p_{X,x}$, $X^d$ equidimensional, choose
$d+p-1$ elements which generate a reduction $K$ of $M$. Denote the matrix 
whose columns are the $d+p-1$ elements by $[K]$; $[K]$ induces
a section of $\Hom \hskip 2pt(\IC^{d+p-1},\IC^p)$ which is a trivial 
bundle over $X$. Stratify $\Hom \hskip 2pt(\IC^{d+p-1},\IC^p)$ by rank. 
Let $[\epsilon]$ denote a $p\times (d+p-1)$ matrix, whose entries are 
small, generic constants. Then, on a suitable neighborhood $U$ of $x$ 
the section of 
$\Hom \hskip 2pt(\IC^{d+p-1},\IC^p)$ induced from $[K]+[\epsilon]$ has at 
most kernel rank 1, is transverse to the rank stratification, 
and the number of points where the kernel rank is 1 is $e(M)$.

\pf The first step is to explain by construction what we mean by ``generic 
constants". Consider the family of maps $G_a$ from 
$X^d$, parametrised by $\IC^{p(d+p-1)}$ to
$\Hom \hskip 2pt(\IC^{d+p-1},\IC^p)$ defined by 
$G_a(x)=G(x,a)=[K(x)]+[A]$, where $[A]$ is the $p\times (d+p-1)$ matrix 
whose entries are coordinates
$a_{i,j}$ on $\IC^{p(d+p-1)}$. Let $\tilde X$ be a resolution of $X$, so 
we have an induced family of maps $\tilde G$ on $\tilde X$. 
Since the map $\tilde G(x,a)$ is a submersion, it follows that for a 
Z-open subset $V$ of $\IC^{p(d+p-1)}$, that for $a\in V$, the map
$\tilde G_a$ is transverse to the rank stratification. We claim that the 
points of $V$ are the generic constants in the theorem.
 Note that the points of $\Hom \hskip 2pt(\IC^{d+p-1},\IC^p)$ of kernel 
rank 1 have codimension $1\cdot((d+p-1)-(p-1))=d$; so since 
$\tilde G_a$ is transverse it can only hit points of the rank 
stratification of kernel rank 1, and only if $D\tilde G_a$ has maximal 
rank at such points which implies
$X$ is smooth at the projection of such points. Let $\tilde K$ be the 
submodule of ${\cal O}^p_{X\times \IC^{p(d+p-1)}}$ defined by the matrix
$[K(x)]+[A]$.

Now apply the multiplicity-polar theorem to $X\times \IC^{p(d+p-1)}$, 
thought of as a family parametrised by $\IC^{p(d+p-1)}$, and 
$(\tilde K,\O^p_{X\times \IC^{p(d+p-1)}})$. Use a point of $V$ as the 
generic parameter value $\epsilon$. 
Then $\O^p_{X\times \IC^{p(d+p-1)}}$ has no polar, because it is free, 
$\tilde K$ has no polar,
 because
$\tilde K$ is generated by $d+p-1$ elements. Choose $U$ a neighborhood of 
$x\times \IC^{p(d+p-1)}$ sufficently small such that every component
 of the cosupport of $\tilde K$ which meets $U$ has $(x,0)$ in its 
closure. Now at $\epsilon$ the cosupport of $\tilde K_{\epsilon}$ is just 
the points where
  $[K]+[\epsilon]$ has less than maximal rank. At such points $e(\tilde 
K_{\epsilon})$ is $1$, because since we are at a smooth point of $X$, the 
local ring of $X$ is Cohen-Macaulay,
so $e(\tilde K_{\epsilon})$ is just the colength, which is 1. Hence 
$e(M)=e(K)=e(\tilde K_{0})=e(\tilde K_{\epsilon})$, which is the number of 
points 
 where the kernel rank of $[K]+[\epsilon]$ is 1.

\rmk 3 In \cite{F} p254 Fulton describes the k-th degeneracy class 
associated to $\sigma$ a homomorphism of vector bundles over $X^d$. The 
support of the class is the set of points where
the rank of $\sigma$ is less than or equal to $k$. Suppose $\sigma: E\to 
F$ where the rank of $E$ is $e$ and the rank of $F$ is $f$, $e\ge f$, 
$e-f+1=d$. Then the $f-1$ degeneracy class is 
supported at isolated points. Fulton shows that if $X$ is Cohen-Macaulay 
at $x$, the contribution to the class at $x$ is the colength of the ideal 
of maximal minors of the
matrix of $\sigma$ at $x$  for some suitable local trivializations of $E$ 
and $F$. Note that this is just the Buchsbaum-Rim multiplicity of the 
module generated by the columns
of the matrix  associated to $\sigma$. Theorem 1.2 shows that in this 
situation if $X$ is pure dimensional, 
the contribution to the degeneracy locus is always the Buchsbaum-Rim 
multiplicity associated to $\sigma$ at $x$, the Cohen-Macaulay hypothesis 
is unnecessary. (Just use the proof of 1.2 to construct a rational 
equivalence to go back to Fulton's case 
close to $x$.)

\sctn  Hypersurface singularities with 1-dimensional singular locus.

In his thesis (\cite{P}) Pellikaan studied non-isolated hypersurface 
singularities. 
This is the setup for his work. He assumed that $f:\IC^{n+1}\to\IC$, $f$ 
had a 1-dimensional singular locus $\Sigma$, which is a complete 
intersection curve defined by an ideal $I$. 
He assumed that $f\in I^2$. This ensured that $J(f)$, the jacobian ideal 
of $f$ was in $I$ as well. (In fact for the singular locus a complete 
intersection Pellikaan proved that if $f$ and its partials were in 
$I$ then $f$ was in $I^2$.) One of the key invariants of $f$ was 

$$j(f)=\dim_{\IC}{{I}\over{J(f)}}$$

\noindent which plays the same role in Pellikaan's work as the dimension 
of ${{\O_{n+1}}\over{J(f)}}$ does in the case of isolated singularities.

Two important examples of  non-isolated singularities are germs of type 
A$_{\infty}$ which have the normal form 
$f(z_1,\dots,z_{n+1})=\sum\limits_{i=1}^{n} {\mathop z_i^2}$ and germs of 
type
D$_{\infty}$ which have normal form 
$f(z_1,\dots,z_{n+1})=z_1z^2_2+\sum\limits_{i=3}^{n+1} {\mathop z_i^2}$. 
Note that if n=2 then D$_{\infty}$ is just a Whitney umbrella.
For A$_{\infty}$ germs $j(f)=0$ while for D$_{\infty}$ germs $j(f)=1$.

Using these building blocks, Pellikaan was able to give a nice geometric 
description of $j(f)$.

\thm 1 Suppose $f$ is as above and $j(f)$ finite. Then $f$ has a 
deformation $F$ such that $F_y$ has $\Sigma_y$ as singular locus for 
generic $y$ where $\Sigma_y$ 
is the Milnor fiber
of $\Sigma$, with only A$_1$ singularities off $\Sigma_y$ and only 
A$_{\infty}$ singularities at points of $\Sigma_y$, except for isolated 
D$_{\infty}$ points. Moreover
$$j(f)=\#\{{\rm D}_{\infty}(F_y)\}+\#\{{\rm A}_{1}(F_y)\}.$$

\pf Cf. \cite {P} p 87 proposition 7.20.

In applying the theory of integral closure to ambient stratification 
conditions like \Af or \Wf in Pellikaan's situation, we see that there are 
three strata--the open stratum, $\Sigma-0$ and the origin. So, there are 
two pairs of ideals $(I,\O_{n+1})$ and $(J(f),I)$ that we are interested 
in. We wish to give 
a geometric interpretation of $(J(f),I)$ using Pellikaan's theorem and the 
multiplicity-polar theorem. First we look at our building block germs.

\prp 2 If $f$ is a germ of type A$_{\infty}$, then $e(I,J(f))=0$, if $f$ 
is a germ of type D$_{\infty}$, then $e(I,J(f))=1$.

\pf If $f$ is a germ of type A$_{\infty}$, then $I=J(f)$, so 
$e(J(f),I)=0$. So suppose $f$ is a germ of type D$_{\infty}$. We may 
assume $f$ is in normal form, as changes of coordinates
do not affect the multiplicity of the pair. We have to compute a sum of 
intersection 
numbers:

$$e(J(f),I)=\Sum_{j=0}^{n} \int D_{J(f),I} \cdot l^{n-j}_{J(f)}\cdot 
l^j_{I}.$$

Consider the part of the sum of form:

$$\Sum_{j=1}^{n} \int D_{J(f),I} \cdot l^{n-j}_{J(f)}\cdot 
l^j_{I}=\Sum_{j=0}^{n-1} \int (D_{J(f),I}\cdot l_{I})\cdot 
l^{n-1-j}_{J(f)}\cdot l^j_{I}.$$  

This is $e(J(f),I)$ where both ideals are restricted to the codimension 1 
polar variety of $I$. Consider the family of candidate polar varieties 
defined by
$z_2=\sum\limits_{i=3}^{n+1} {a_iz_i}$. Since this a Z-open subset of all 
potential polar varieties, if we show that for a Z-open subset of them 
that the multiplicity of the 
pair of the restriction of the 
ideals to each candidate in the set is zero then we will have shown that 
all of terms in this second sum are zero and all these candidates are 
actually polars. 
Now it is obvious from the normal form of $f$ that when we restrict our 
two ideals to any element of this set the two ideals become equal so all 
of the terms in the second sum are zero.

It remains to compute $ \int D_{J(f),I} \cdot l^{n}_{J(f)}$. Our approach 
is to choose a Z-open set of candidate polar curves of $J(f)$, 
then show that each candidate gives the same value for the computation of 
the desired intersection number. Consider the family of curves defined by 
ideals
$J_{a,b,c}=(b_1(z_1z_2)+c_1z^2_2+\sum\limits_{i=3}^{n+1} 
{a_{1,i}z_i},\dots, b_n(z_1z_2)+c_nz^2_2+\sum\limits_{i=3}^{n+1} 
{a_{n,i}z_i}) $. 
For a Z-open set of coefficients, we can re-write the ideals defining 
these curves as 
$$J_{a,b,c}=(z_1z_2+cz^2_2,\dots, z_i+b_iz_1z_2,\dots) $$
\noindent where $3\le i\le n+1, c\ne 0$.
Each curve in this family has two components; one of which (given by 
$z_2=0$) lies in $V(J(f))$. The other component is the candidate polar 
curve. So we get the family of
parmeterizations
$\phi(t)=( -ct,t,\dots, b_ict^2,\dots)$ for the candidate polar curves. 
Now the intersection number we want is just the multiplicity of the pair 
restricted to a polar curve; by the additivity
of the multiplicity (\cite{K-T}) this is just $e(J(f))-e(I)$ restricted to 
the polar curve; given a parameterization this is just the order of 
vanishing of $\phi^*(J(f))$ less
the order of vanishing of $\phi^*(I)$. Now $\phi^*(J(f))=(t^2)$ and 
$\phi^*(I)=(t)$ for all parameterizations, so the value of this 
intersection number is $2-1=1$, so $e(J(f),I)=1$.

For our basic building block germs we have seen that $j(f)=e(J(f),I)$. The 
next theorem shows that this is true in general. If $F$ depends on 
coordinates $(y,z)$, let $J_z(F)$ denote the ideal generated by the 
partials of $F$ with respect to $z$.

\thm 3 Suppose $f:\IC^{n+1}\to\IC$, $f$ has a 1-dimensional singular locus 
$\Sigma$, which is a complete intersection curve defined by an ideal $I$, 
$f\in I^2$ and $j(f)$ finite.
Then $$j(f)=\dim_{\IC}{{I}\over{J(f)}}=e(J(f),I).$$

\pf Let $F$ be the deformation of Theorem 2.1. Denote the parameter space 
by $Y^k$. The singular set of $F$ is given by a complete intersection 
$\tilde I$. We are interested 
in the family of pairs of ideals given by $(J_z(F),\tilde I)$ as these 
restrict to $(J(f) ,I)$ at $y=0$. Since $\tilde I$ defines a complete 
intersection it
 has no polar variety of dimension $k$. 
Since $J_z(F)$ is generated by $n+1$ generators it has no polar of 
dimension $k$ either. This means that the multiplicity of the pair at the 
origin is same as 
the sum of the multiplicities over a generic parameter value by the 
multiplicity-polar formula. Pick a generic $y$. We have 
$(J_z(F))_y=J(F_y)$, so the cosupport of $(J_z(F))_y$
consists of A$_1$ points off $\Sigma_y$, isolated D$_{\infty}$ points on 
$\Sigma_y$ and A$_{\infty}$ points. Off $\Sigma_y$, $(\tilde 
I)_y=\O_{n+1}$, so off
 $\Sigma_y$, at A$_1$ points,  $e(J(F_y),\tilde I_y)=e(J(F_y),\O_{n+1})=1$ 
and $0$ elsewhere off $\Sigma_y$. On $\Sigma_y$, $e(J(F_y),\tilde I_y)=1$ 
at D$_{\infty}$ points,
otherwise it is $0$ by proposition 2.2. So the sum of the $e(J(F_y),\tilde 
I_y),z)$ at points where it is non-zero is just $\#\{{\rm 
D}_{\infty}(F_y)\}+\#\{{\rm A}_{1}(F_y)\}$.

Then, by the multiplicity-polar formula we know that $$e(J(f),I)=\#\{{\rm 
D}_{\infty}(F_y)\}+\#\{{\rm A}_{1}(F_y)\}$$
which proves the theorem.

If $R$ is Cohen-Macaulay of dimension d, $M$ a submodule of a rank $p$ 
free module $F$ of finite colength, then by a theorem of Buchsbaum and Rim 
(\cite{B-R}), $e(M,F)$, which is $e(M)$, is just the colength
of $M$ if $M$ has $d+p-1$ generators. Theorem 3 can be viewed as a first 
step in generalizing this result to pairs of modules.

Using some other results of Pellikaan, we can link $e(J(f),I)$ and the 
L\^e numbers introduced by Massey. In the situation of Theorem 2.3 there 
are two L\^e numbers--$\lambda^0(f)$ and $\lambda^1(f)$; denote the number 
of D$_{\infty}$ points of $f$ by $\delta(f)$.

\prp 4 Assume the hypotheses of Theorem 2.3, then 
$$\lambda^0(f)=e(J(f),I)+e(JM(\Sigma))+\delta(f)$$

\pf If $F$ is the Milnor fiber of $f$, we have that 
$$\chi(F)=1+(-1)^{n-1}\lambda^1(f)+ (-1)^n\lambda^0(f)=1+(-1)^n(j(f)+ 
\delta_f+\mu(\Sigma)-1).$$

The first equality is due to Massey (\cite{Ma}), while the second is due to Pellikaan 
(\cite {P}, p113, proposition 10.11).
In the present situation, since the transverse Milnor number is 1, 
$\lambda^1(f)={\rm mult}\hskip 2pt(\Sigma)$, 
while $e(JM(\Sigma))=\mu(\Sigma)+{\rm mult}\hskip 2pt(\Sigma)-1$.

Therefore, substituting and canceling we get
$$\lambda^0(f)=e(J(f),I)+e(JM(\Sigma))+\delta(f).$$

Now we turn to the extension of these ideas to hypersurface singularities 
with a higher dimensional singular locus.

\sctn Hypersurface singularities with d-dimensional singular locus.

In this section we assume that $I=(g_1,\dots,g_p)\subset\O_n$ defines a 
complete intersection of dimension $d>1$, and $S(f)=V(I)$, hence we can 
write $f$ as
$f=\sum_{i,j}^ph_{i,j}g_ig_j$, where $h_{i,j}=h_{j,i}$, for some 
$h_{i,j}$. Let $[H]$ denote the symmetric matrix with entries $h_{i,j}$. 
We will want to study those germs $f$
for which
$j(f)<\infty$. Basic examples of such germs are those of type $A(d)$. For 
these germs up to a change of coordinates, $I=(z_1,\dots, z_{n-d})$, 
$f=\sum_{i=1}^{n-d}z_i^2$,
$z_i$ part of a coordinate system on $\IC^n$. It turns out that the 
condition that $j(f)<\infty$ is much more restrictive than in the case 
where dimension of
$V(I)=1$. Pellikaan already showed that $j(f)<\infty$ implies $I$ defines 
an ICIS. The next proposition gives a further restriction.

\prp 1 Suppose $f$, $I$ as above, then if $[H]$ has less than maximal rank 
at the origin, the set of points on $V(I)$ where the singularity type is 
not $A(d)$ is of
codimension 1 in $V(I)$ , hence $j(f)$ is not finite.

\pf If $f$ has an $A(d)$ singularity at $x\in V(I)$ then $V(I)$ is smooth 
at $x$ and the matrix $[H(x)]$ must have rank $n-d$. But the points where 
$\det [H]=0$ defines a non-empty hypersurface in $V(I)$, since $\det 
[H(0)]=0$ and the dimension of $V(I)>1$. Hence,
 at these points  $f$ does not have an $A(d)$ singularity. Since at these 
points $I\ne J(f)$, it follows that $j(f)=\infty$.

There are two types of L\^e cycles; those which are the images in $\IC^n$ 
of components of the exceptional divisor of the jacobian blow up, called 
fixed cycles, and the polar
varieties of the fixed cycles called moving cycles.

\cor 2 Suppose $f$, $I$, $[H]$ as above. Then $V(I)$ contains a fixed L\^e 
cycle of dimension $d-1$.

\pf  We can deform $f$ so that the $D^{\infty}$ points are dense in the zero set of $\det [H]=0$. These points  are clearly the image of a component of the exceptional divisor by Proposition 2.4, and by the properties of the L\^e numbers. Then when we specialize, the component of $E$ will specialize as well.

In \cite {P} and in \cite{P1} Pellikaan defines the singularities of type 
$D(d,p)$; here $d$ is the dimension of $S(f)$, 
while $p$ is the dimension of the kernel rank of $[H]$
at the point in question. Then  $f:\IC^n,x\to \IC,0$ has type $D(d,p)$ at  
$x$ if local coordinates can be chosen so that $f$ has the local form
$$f=z_1^2+\dots+z_q^2+\sum_{1\le i\le j\le p} x_{i,j}y_iy_j$$
\noindent where $z,x,y$ are part of a coordinate system on $\IC^n$ at $x$, 
$n-d=q+p$. From 3.1 it follows that if $f$ has singularity type $D(d,p)$ 
at the origin, and $d>1$,
then, since $\det[H(0)]=0$, it follows that $j(f)=\infty$, contrary to 
remark 5.3 of \cite {P} and Remark 5.4 of \cite {P1}. This  shows that 
$j(f)$ fails to be finite
in what seems to be the next most simple case to the $A(d)$ singularities 
when $d>1$. Instead, the structure of $S(f)$ seems more like a 
discriminant, in that the non-generic
points appear in codimension 1.

In the next lemma we begin to look at those germs where $[H]$ has maximal 
rank, so we can characterize those germs where  $j(f)<\infty$.

\lem 3 Suppose $f=\sum_{i,j}^ph_{i,j}g_ig_j$, $\det [H(0)]\ne 
0$, $I=(g_1,\dots,g_p)\subset\O_n$. Then one can chose a set of
generators
$(g'_1,\dots,g'_p)$ of $I$ such that $f=\sum_i^p (g'_i)^2$.

\pf The proof is standard, so we just sketch the details. Given an 
invertible matrix $[R]$ with entries in $\O_n$, it is clear that if 
$$[g]=[R][g'],$$
\noindent where $[g]$ is the column vector whose entries are the $g_i$, 
$[g']$ another column vector, that the entries of $[g']$ are also a set of 
generators of $I$. Given

$$[f]=[g]^t[H][g]$$

and 

$$[g]=[R][g']$$

\noindent it follows that 

$$[f]=[g']^t([R]^t[H][R])[g'].$$

Hence, we need to show that by choice of $[R]$ we can reduce $[H]$ to the 
identity matrix.
This is done in two steps--first we can chose $[R]\in Gl(p,\IC)$ so that 
we can assume $[H(0)]=I$. (This follows because the action of $Gl(p,\IC)$ 
clearly preserves rank,
 the orbits of $Gl(p,\IC)$ are connected constructible sets, and the 
orbits of non-singular matrices are open, by a tangent space calculation.)

For the second step we assume $[H(0)]=I$, consider the linear homotopy 
from $I$ to $[H]$; this stays inside the set of invertible symmetric 
matrices. The congruence
transformation gives an action of the group ${\cal C}$ of invertible 
$p\times p$ matrices with entries in $\O_n$ on the $p\times p$ symmetric 
matrices. Applying the
techniques of Mather-Damon produces a homotopy in ${\cal C}$ which 
trivializes our linear homotopy, which finishes the proof.

Lemma 3.3 also appears as a remark without proof in \cite{Z} (see page 87).

Given a set of generators $\{g_1,\dots,g_p\}$ for an ideal, we can form 
the function $G$ whose components are the $g_i$. If $\{g_1,\dots,g_p\}$ define an ICIS, then the map $G$ is said to be of {\it finite singularity type}.

\cor 4 Suppose $f=\sum_{i,j}^ph_{i,j}g_ig_j$, $\det [H(0)]\ne 
0$,$I=(g_1,\dots,g_p)\subset\O_n$, $\{z_i\}$ coordinates on $\IC^p$ 
then generators $(g'_1,\dots,g'_p)$ of $I$ can be
chosen so that 
$$f=\sum_{i=1}^p z_i^2\circ G'$$

\pf By lemma 3.3 we have there exists generators $(g'_1,\dots,g'_p)$ of 
$I$ such that 
$$f=\sum_{i=1}^p (g'_i)^2=\sum_{i=1}^p z_i^2\circ G'.$$

Thus the study of functions with $j(f)<\infty$ is intimately tied up with the study of 
functions on the discriminant of a map germ of finite singularity type as we shall see below.

We wish to describe a condition which will ensure that the pullback by $G$ 
of a function on $\IC^p$ with a Morse singularity at the origin gives a 
function on $\IC^n$ with
$j(f)<\infty$ for the ideal defined by the components of $G$. This 
completes our geometric description of the meaning of $j(f)$ finite.

Our condition is based on the intersection of the levels of the Morse 
function in the target with the discriminant, $\Delta(G)$, of $G$. At this 
point we asssume that $I$
defines an ICIS. This implies that if $G$ comes from a minimal set of 
generators of $I$, then $G|S(G)$ is a finite map.

We can partition $S(G)$ by the  $S_i(G)$ which denotes points of $S(G)$ 
where the kernel rank of $G$ is $i$. We can also partition $\Delta(G)$ as 
follows.
For each point $z$ of $\Delta(G)$, list the points $S_z$ of $S(G)$ mapped 
to $z$. The points $z$ and $z'$ are in the same element of the partition 
if there is a bijection
between $S_z$ and $S_{z'}$ which which preserves components of the 
$S_i(G)$. It is easy to see that the elements of this partition are 
constructible sets since $G|S(G)$ is
finite. Given an element of the partition of $\Delta(G)$, we now associate 
a collection of systems of linear sub spaces of $T\IC^p$ over the 
underlying set $P$ of the partition
element. Since
$G$ has constant rank on each $S_i(G)$,  $D(G)|_{S_i(G)\cap 
G^{-1}(P)}(T\IC^n|S_i(G)\cap G^{-1}(P))$ is a well defined sub bundle of     
$G^*T\IC^p$ over $S_i(G)\cap G^{-1}(P)$. Since the restriction of $G$ to   
each component of $G^{-1}(P)$ is a homeomorphism or finite
cover, the push forward by $G$ of these sub bundles gives the desired 
collection of systems of linear spaces. We call the partition
of $\Delta (G)$  together with the collection of linear spaces on each 
element of the partition an {\it enriched partition}.  A smooth
subset $V$ of $\IC^p$ is {\it enriched transverse} to the enriched 
partition if at every point of intersection with the elements of the 
partition the
tangent space of $V$ is transverse to each of the linear spaces we have 
associated to the element of the partition at that point.  Since   
the restriction of $G$ to each component of $G^{-1}(P)$ is a homeomorphism 
or finite
cover, all of the linear spaces at a smooth point in a partition element 
contain the tangent space to the partition element. So if $V$ is
transverse to each element of the partition it is enriched transverse. The 
next proposition describes a situation in which transversality and 
enriched transversality are
equivalent.

\prp 5 Suppose there exists an element $P$ of the partition which is a 
Z-open subset of $\Delta (G)$ whose pre-images lie in the Z-open subset 
$S_{n-p+1}(G)$ on which
$G$ is immersive.  Then all of the systems of linear spaces associated to 
$P$ are just the tangent bundle to $P$. 

\pf Suppose  $y\in P$, $z$ a preimage in  $S_{n-p+1}(G)$. Since $G$ 
restricted to $S_{n-p+1}(G)$ is immersive at $z$, the dimension of 
$DG(TS_{n-p+1}(G))$ is $p-1$ which is
the dimension of $D(G)(z)T\IC^n$, so these spaces are equal; further 
$DG(TS_{n-p+1}(G))$ is the tangent space to $\Delta (G)$ at $y$, which is 
the tangent space to $P$ at
$y$.

Now we give our condition for $j(f)$ finite.

\thm 6 Suppose  $I=(g_1,\dots,g_p)\subset\O_n$ defines an ICIS of 
dimension $d>1$, $G$ the mapgerm whose components are the $g_i$, 
$h:\IC^p,0\to \IC,0$ a function with an
isolated singularity at the origin, $f=h\circ G$. Then
$$j(f,G^*(J(h))\O_n):={\rm dim}_{\IC}{{G^*(J(h))\O_n}\over{J(f)}}<\infty 
$$ \noindent
if and only if $h^{-1}(0)$ is enriched transverse to the enriched 
partition of $\Delta(G)$ except possibly at the origin.

\pf Suppose $j(f,G^*(J(h))\O_n)$ finite. Then, except possibly at the 
origin, $J(f)=G^*(J(h))\O_n$. If the enriched transversality condition 
fails, there must be a curve
$\phi:\IC\to \Delta G$, such that the image of $\phi$ lies in an element 
of the partition, and the tangent space to $h^{-1}(0)$ contains one of the 
systems of linear spaces
along the partition element. This implies that contained system is in the 
kernel of $Dh$ along $\phi$. Then $\phi$ has a lift to the component of 
$S_i(G)$ associated to the
contained system, denoted
$\psi$. Along the image of $\psi$ we have

$$Df\circ \psi(T\IC^n)=Dh\circ(G\circ \psi)DG\circ\psi(T\IC^n)=0$$

Hence $V(J(f))\supset {\rm im}\psi$, while $V(G^*(J(h))\O_n)=V(I)$ which 
is a contradiction.

Suppose enriched transversality holds. If  $j(f,G^*(J(h))\O_n)$ is not 
finite, there exists a curve $\psi$ whose image properly contains the 
origin in $\IC^n$, such that
$J(f)\ne G^*(J(h))\O_n$ along $\psi$. At points of $\IC^n$ off $S(G)$, $G$ 
is a submersion, hence $J(f)=G^*(J(h))\O_n$. If $\psi$ lies in $S(G)$, 
then the image of $\psi$
lies in $S(f)$ since $G^*(J(h))\O_n=\O_n$ at such points. Then $\psi$ lies 
in the zero set of $F$, hence $G\circ \psi$ lies in the zero set of $h$. 
Then enriched transversality
 fails along $G\circ \psi$.

\cor 7  Suppose $h$ has a Morse singularity at the origin in the set-up of 
Theorem 3.6, then $j(f)$ is finite if and only if  $h^{-1}(0)$ is enriched 
transverse to the
enriched partition of $\Delta(G)$ except possibly at the origin.

\pf If $h$ has a Morse singularity, then $ G^*(J(h))\O_n=I$.

\cor 8 Suppose $I=(g_1,g_2)$ in the setup of Theorem 3.6. Then 
$j(f,G^*(J(h))\O_n)$ is finite iff $f^{-1}(0)\cap S(G)$ is the origin.

\pf If $p=2$, then $\Delta (G)$ is a curve, and $G$ restricted to each 
branch of $S(G)$ is an immersion except at the origin. Then enriched 
transversality becomes ordinary
transversality, so $h^{-1}(0)$ must miss $\Delta (G)$ off the origin, so 
$f^{-1}(0)\cap S(G)$ is the origin.

Theorem 3.6 introduces an interesting class of functions. Given an ICIS, 
by using appropriate $h$ we can construct examples of non-isolated 
singularities in which the
singular locus is the ICIS, but the transverse singularity type is 
constant and is that of $h$. In studying the equisingularity of families 
of such examples, the key invariant
is the multiplicity of the pair $J(f),G^*(J(h))\O_n$. This number should 
also be linked to the way $h^{-1}(0)$ meets the discriminant of $G$ at the 
origin.

 Now we show that such functions with $j(f)$ finite 
are plentiful.

\prp 9 Suppose $G:\IC^n,0\to\IC^p,0$, $G^{-1}(0)$ an ICIS, $p>1$. Then if 
$h_a(x)=\sum a_ix^2_i $, for $a\in U$, $U$ a Z-open subset of $\IC^p$, 
$h^{-1}(0)$ is transverse to
the enriched partition of $\Delta (G)$ except perhaps at the origin.

\pf Consider $H(a,z)=\sum a_ix^2_i\circ G(z)$. We have 

$$DH=<\dots,x^2_i\circ G,\dots, 2a_ix_i\circ G,\dots>$$
This implies that $H$ is a submersion except along $\IC^p\times 
G^{-1}(0)$. Denote by $\pi$ the projection of $H^{-1}(0)$ to $\IC^p$. By 
Sard's lemma for varieties (prop. 3.7
p42 \cite{M}) there exists a Z-open subset $U\subset \IC^p$ such that $\pi$ 
is smooth at $z\in H^{-1}(0)\cap\pi^{-1}(U)/\IC^p\times G^{-1}(0)$. This 
implies that the fiber of
$\pi$, which is the fiber of $h_a\circ G$ over $0$ is smooth at $z$; in 
addition since $\pi$ maps $T_z(H^{-1}(0))={\rm ker}DH_z$ surjectively to 
$\IC^p$, the ${\rm ker}DH_z$
does not contain $\IC^n$, thus $h_a\circ G$ is a submersion at $z$ as 
well, hence $h_a$ is enriched transverse to the enriched partition of 
$\Delta(G)$, except perhaps at the
origin.

Now that we know that it is worth proving results about functions with 
$j(f)$ finite for $V(I)$ an ICIS of dimension $>1$, we prove the analogue of 
2.3. To do this we first
study a special deformation of $f=\sum_1^p z^2_i$.

We call the following pair of deformations a smoothing of $f$.

$$F(u,b,z)=\sum_{i}(1+\sum_j b_{i,j}z_j) (g_i-u_i)^2$$

$$\tilde G(u,z)=(g_1(z)-u_1,\dots,g_p(z)-u_p)$$

This is called a smoothing because of the following lemma:

\lem 10 For a Z-open subset $U$ of $\IC^p\times \IC^{pn}$, $f_{u,b}$ has 
only A$_1$ singularities off $G^{-1}(u)$, $G^{-1}(u)$ is smooth and 
$f_{u,b}$ has only $A(d)$
singularities on $G^{-1}(u)$.

\pf let $V\subset \IC^p$ be the complement of $\Delta(G)$ in $\IC^p$, then 
$G^{-1}(u)$ is smooth for $u\in V$.

 We claim $D_zF(u,b,z)$ is a submersion off $\IC^{pn}\times\Gamma(G)$, 
where
$\Gamma(G)\subset\IC^n\times\IC^p$ denotes the graph of $G$.

Let $e_i$, where $1\le i\le n$ denote the unit vectors in $\IC^n$. Then we 
have 

$$\Pd{D_z(F)}{b_{i,j}}=(g_i-u_i)^2e_j.$$

This implies $D_zF(u,b,z)$ has maximal rank when some $g_i-u_i$ is not 
zero which proves the claim.

Now consider $D_zF(u,b,z)^{-1}(0)$. The claim shows this is smooth off  
$\IC^{pn}\times\Gamma(G)$.

As in the proof of 3.9 we consider the projection of this set to 
$\IC^p\times \IC^{pn}$, let $W$ be the Z-open subset of the base over 
which $\pi$
is smooth off $\IC^{pn}\times\Gamma(G)$. Now the tangent space to 
$D_zF(u,b,z)^{-1}(0)$ at a point $x$ is just the kernel at $x$ of 
$D(D_zF(u,b,z))$, which has dimension
$p+pn$ and which surjects to $\IC^p\times \IC^{pn}$. Hence $D^2_zF(u,b,z)$ 
has maximal rank, so $f_{u,b}$ has only Morse singularities off $G=u$. 
 Let $U=W\cap \IC^{pn}\times V$, then for $(u,b)\in U$, we have $g_u$ has 
a smooth fiber over zero. Since the set of points where the matrix $H$ 
with entries 
$h_{i,i}=\sum_j 1+b_{i,j}z_j$, $h_{i,j}=0$ $i\ne j$ has maximal rank on 
some Z-open subset of $\IC^{pn}\times\IC^n\times\IC^p$ which contains 
zero, we can ensure that each of 
the $f_{u,b}$ has only $A(d)$ singularities on some fixed neighborhood of 
the origin in $\IC^n$ on $g_u=0$.

{\it Remark}  It was pointed out to the author by the referee that this lemma also follows from the statement and proof of theorem 1 of 
Now we extend Theorem 2.3 to ICIS of dimension greater than 1.

\thm 11 Suppose $f:\IC^{n+1}\to\IC$, $f$ has a d-dimensional singular 
locus $\Sigma$, $d>1$, which is an ICIS defined by an ideal $I$, $f\in I^2$ 
and $j(f)$
finite. Then $$j(f)=\dim_{\IC}{{I}\over{J(f)}}=e(J(f),I)=\#A_1(f).$$
\noindent where $\#A_1(f)$ is the number of $A_1$ singularities appearing 
in a smoothing of $f$.

\pf The proof is similar to that of 2.3. By \cite{P} Theorem 3.1 p145 the 
quotient of the ideals $(g_1(z)-u_1,\dots,g_p(z)-u_p)/J_z(F)$ is perfect, 
where F is part of a
smoothing of $f$, hence the length of the quotients $(g_u)/J(f_{u,b})$ is 
independent of parameter, and for generic parameter value is just 
$\#A_1(f)$. Meanwhile, $J_z(F)$ and 
$(g_1(z)-u_1,\dots,g_p(z)-u_p)$ have no polar varieties of dimension 
$p+p(n+1)$, so as in Theorem 2.3, the multiplicity polar theorem implies 
that $e(J(f_{u,b},(g_u))$ is
independent of parameter, so again is $\#A_1(f)$, hence the theorem 
follows.

Now we wish to extend proposition 2.4 to ICIS of dimension greater than 1.

In \cite {Z}, prop 5.5.5, p86, (cf. also \cite {N}), Zaharia computed the homology of the Milnor 
fiber, $\hat f$, of a function germ $f$ defined on
$\IC^{n+1}$ whose singular set $\Sigma$ was an ICIS of codimension $p$ 
such that $j(f)<\infty$. His result was:

$$H_*(\hat f)=\left\{ {\matrix{Z, \hskip 2pt{\rm if}\hskip 2pt *=0,p-1 \cr
Z^{\mu_{\Sigma}+\sigma} \hskip 2pt{\rm if}\hskip 2pt *=n\cr
0, \hskip 2pt {\rm otherwise}\cr
}} \right.$$

Here $\sigma$ is the number of $A_1$ points appearing in a smoothing, 
which we have shown is $j(f)$.

\prp 12 Assume the hypotheses of Theorem 3.11, then 
$$\lambda^0(f)=e(J(f),I)+e(JM(\Sigma))$$

\pf By Massey (\cite{Ma}) we have that 

$$\chi(\hat f)=1+\sum_{i=0}^d (-1)^{n-i}\lambda^i(f)$$

$$=1+(-1)^n\lambda^0(f)+\sum_{i=1}^d(-1)^{n-i}\lambda^i(f)$$
Now, for $i>0$, since $\Sigma$ is the only fixed L\^e cycle of dimension greater than $0$, and $f$ has transverse Milnor number $1$, since the type of $f$ is $A(d)$
generically on $\Sigma$, 

$$\lambda^i(f)=m_{d-i}(\Sigma),$$

\noindent where $m_{d-i}(\Sigma)$ is the $d-i$ polar multiplicity of the ICIS $\Sigma$. In turn, $m_{d-i}(\Sigma)=\mu^{d-i}(\Sigma)+\mu^{d-i-i}(\Sigma)$ (\cite{G-4})
where $\mu^{d-i}(\Sigma)$ is the Milnor number of $\Sigma\cap H_i$ where $H_i$ is a generic plane of codimension $i$, and where $\mu^{-1}=1$.

Substituting, the sum telescopes to: 

$$\chi(\hat f)=1+(-1)^n\lambda^0(f)+(-1)^{n-d}+(-1)^{n-1}\mu^{d-1}(\Sigma).$$

Calculating $\chi(\hat f)$ from the homology calculation of \cite{Z} we get:

$$1+(-1)^n\lambda^0(f)+(-1)^{n-d}+(-1)^{n-1}\mu^{d-1}(\Sigma)=1+(-1)^{n-d}+(-1)^n(\mu_{\Sigma}+\sigma),$$
Hence

$$\lambda^0(f)=\sigma+\mu_{\Sigma}+\mu(\Sigma\cap H_1)$$

$$=e(J(f),I)+e(JM(\Sigma)).$$

\rmk 13 There are  two other general calculations of the homolgy of the Milnor fiber  in $\cite {Z}$ (Theorem 5.5.4 
and Proposition 5.5.6). (Note, however the typo in the formula of 5.5.4--the
coefficients of $\mu_{\Delta}$ and $\mu_{\Sigma}$ should be exchanged.) Using these calculations, it is possible to prove by the same methods as 3.12, two other formulas for
$\lambda^0(f)$. 

In the first case, assume $V(I)=\Sigma$ is an ICIS of dimension 2, write $[f]=[g]^t[H][g]$ as we did earlier, let $H$ denote the ideal generated by $I^2$ and the entries
of $[H][g]$, assume $\dim_{\IC}H/J(f)$ is finite, $V(\det [H])\cap\Sigma=\Delta$, where $\Delta$ is an ICIS of dimension 1. We can consider
the smoothing used by Zaharia to study
this situation, and the ideal $H$ extends to $\tilde H$ in a natural way, to the space of the smoothing.  Then the polar of $\tilde H$ may be non-empty if the kernel rank of
$[H]$ is $>2$. Call the multiplicity of the polar of $\tilde H$ over the base $m( \Gamma(\tilde H))$. Then the multiplicity polar theorem applied to the smoothing gives
$e(J(f),H)+m( \Gamma(\tilde H))=\#(A_1(f))$ and hence,
$$\lambda^0(f)=e(J(f),H)+m( \Gamma(\tilde H))+e(JM(\Sigma))+2e(JM(\Delta)).$$

In the second case, assume $V(I)$ has dimension $d>1$, assume the rank of $[H(0)]$ is $p-1$ (one less than maximal). Then, as Zaharia remarks (\cite{Z} p. 87), generators
$(g_1,\dots,g_p)$ for $I$ can be found so that $f=\det([H])g^2_1+g^2_2+\dots +g^2_p$. Then the ideal $H$ of the last paragraph is just $(\det([H])g_1, g^2_1,g_2,\dots,g_p)$.
Since $H$ has only $p+1$ generators as does $\tilde H$ the polar of $\tilde H$ is empty and 
$$\lambda^0(f)=e(J(f),H)+e(JM(\Sigma))+2e(JM(\Delta)).$$

The form of these formulae makes it likely that they are special cases of a more general theorem.

It has  long been known that in cases like those considered here, that the independence  from parameter of the L\^e numbers implies that the families $\Sigma(t)$ and
$\Delta(t)$ are Whitney equisingular (See for example \cite {G-G} prop 4.6, for the case where $I=J(f)$, and use the fact that the components of the exceptional divisor of
the blowup of $\IC^{n+1}$ by
$J(f)$ which project to $\Sigma$ and $\Delta$ are the conormals of $\Sigma$ and $\Delta$. ) Thus, a relation between the L\^e numbers and the invariants used to control the
Whitney equisingularity of $\Sigma$ and $\Delta$ is not unexpected. That the formulae relate $\lambda_0$ so simply to the zero dimensional invariants of the strata and to the
\Af invariant is surprising.

Now we develop some results which shows how well
$e(J(f),I)$ is linked to the \Af and \Wf conditions.

\sctn Conditions \Af and \Wf

In this section, we'll study Thom's
Condition A$_f$, and Henry, Merle and Sabbah's Condition \Wf, which
concern limiting tangent hyperplanes at a singular point of a complex
analytic space.  First we recall the notions of integral dependence and 
strict dependence.

Let $(X,0)$ be the germ of a complex analytic space, and $\E:=\O_X^p$
a free module of rank $p$ at least 1. Let $M$ be a coherent submodule
of $\E$, and $h$ a section of $\E$.  Given a map of germs
$\vf\:(\IC,0)\to(X,0)$, denote by $h\o\vf$ the induced section of the
pullback $\vf^*\E$, or $\O_{\IC}^p$, and by $M\o\vf$ the induced
submodule.  Call $h$ {\it integrally dependent\/} (resp., {\it strictly
dependent\/}) on $M$ at $0$ if, for every $\vf$, the section $h\o\vf$
of $\vf^*\E$ is a section of $M\o\vf$ (resp., of $m_1(M\o\vf)$, where
$m_1$ is the maximal ideal of $0$ in $\IC$).  The submodule of $\E$
generated by all such $h$ will be denoted by $\?{M}$, resp., by $\sdep
M$ .

In the context of hypersurface singularities, given a family of map-germs 
$F(y,z)$ parametrised by $Y=\IC^k$, where 
$F:\IC^k\times\IC^{n+1},\IC^k\times 0,0\to \IC,0,0$
Thom's \Af condition holds for the pair
$(\IC^k\times\IC^{n+1}-S(F),\IC^k\times 0)$ at $y\in Y$ if and only if 
every limit of tangent hyperplanes to the fibers of $F$ on 
$\IC^k\times\IC^{n+1}-S(F)$ contains
$TY$ at $y$. The condition holds for the pair if it holds for the pair at 
every $y$. Although this condition  looks like it says nothing about 
strata other than the open stratum,
 this can be deceiving, as we shall see.

\prp 1 Suppose $F:\IC^k\times\IC^{n+1},\IC^k\times 0,0\to \IC,0,0$ then the 
following are equivalent:

1)  The \AF condition holds for the pair 
$(\IC^k\times\IC^{n+1}-S(F),\IC^k\times 0)$ at $0$. 

2) The fiber over $0$ of the exceptional divisor $E$ of the blowup of 
$\IC^k\times\IC^{n+1}$ by $J(F)$, denoted $B_{J(F)}(\IC^k\times\IC^{n+1})$ 
is contained in $C(Y)$, 
the conormal of $Y$.

3) ${{\partial F} \over {\partial y_i}}\in \sdep
{J(F)}$ for $1\le i\le k$.

4)  ${{\partial F} \over {\partial y_i}}\in \sdep
{J_z(F)}$ for $1\le i\le k$.

\pf The fiber over $0$ of the exceptional divisor $E$ of 
$B_{J(F)}(\IC^k\times\IC^{n+1})$ is exactly the set of limiting tangent 
hyperplanes
at $0$
to the fibers of $F$ on $\IC^k\times\IC^{n+1}-S(F)$; saying that this 
fiber lies in the conormal of $Y$ just says that each limit contains 
the tangent space to $Y$ at $0$. This shows 1) and 2) are equivalent. The 
equivalences of 1) and 3) and 4) can be found in \cite{GK}.

The \WF condition holds for the pair 
$(\IC^k\times\IC^{n+1}-S(F),\IC^k\times 0)$ at $0$ if there exist a 
(Euclidean) neighborhood $U$ of $0$
in $\IC^k\times\IC^{n+1}$ and a constant $C>0$ such that, for all $y$ in 
$U\cap Y$ and all
$x$ in $U\cap (\IC^k\times\IC^{n+1}-S(F))$, we have
   $$\dist\bigl\leftp T_yY(F(y)), T_x(\IC^k\times\IC^{n+1})(F(x)\bigr)\leq 
C\,\dist(x,Y)$$
 where $T_yY(F(y))$ and $T_x(\IC^k\times\IC^{n+1})(F(x))$ are the tangent 
spaces to the
indicated fibers of $F$ and the restriction $F|Y$.  

\prp 2 Suppose $F:\IC^k\times\IC^{n+1},\IC^k\times 0,0\to \IC,0,0$ then the 
following are equivalent:

1) The \WF condition holds for the pair 
$(\IC^k\times\IC^{n+1}-S(F),\IC^k\times 0)$ at $0$.

2) ${{\partial F} \over {\partial y_i}}\in \?{\mY J(F)}$ for $1\le i\le 
k$.

3)  ${{\partial F} \over {\partial y_i}}\in \?{\mY J_z(F)}$ for $1\le i\le 
k$.

\pf This follows from proposition 1.1 of \cite{GK2}

Now we want to look at the connection between the multiplicity of the 
pair, $e(J(f),I)$, and the \AF condition.
 At this point we no longer assume that $I$ defines a curve singularity. 
We do need two simple lemmas first.

\lem 3 Suppose $I$ is an ideal generated by $d$ elements in an 
equidimensional local ring $R$ of dimension $n$  such that $R/I$ has 
dimension $n-d$. Suppose
$J\subset I$ is a reduction of $I$. Then $J=I$.

\pf The proof is by induction on $d$. Assume $d=1$, denote the generator 
of $I$ by $p_1$. Let $J=(f_1p_1,\dots,f_kp_1)$. 
 If some $f_i$ is a unit, then we are done.
Suppose no $f_i$ is, and denote the ideal they generate by $F$. If $p_1$ 
satisfies a relation of integral dependence, then 
we get $$(p_1)^k+\sum\limits_{i=0}^{k-1} g_ip^i_1=0$$
\noindent where $g_i\in J^{k-i}$. Then $g_i\in F^{k-i}(p^{k-i})$, so the 
equation of integral dependence implies that 
there exists a unit $u$ such that $u\cdot p^k=0$ which is a contradiction.

Assume $I$ is generated by $d$ elements ; work on $R'=R/(p_1)$, then 
applying the induction hypothesis to the homomorphic images
of $J$ and $I$ in $R'$ we have that these images are equal, hence 
$p_i=g_i+r_ip_1$ where $g_i\in J$. Notice that $\{p_1, p_2-r_2p_1, 
p_3,\dots,p_d\}$ 
is a set of generators for $I$. 
Now mod out by $p_2-r_2p_1=g_i$, and again apply the induction hypothesis. 
This shows that  $\{p_1, p_3,\dots,p_d\}$ are in $J$ hence $I$ is in $J$ 
since the missing generator of
$I$ is already in $J$.

Note that it is not necessary for $I$ to be radical.

We say that $f:\IC^{n+1},x\to \IC,0$ has singularity of type A($d$) at 
$x$, if local coordinates $(z_1,\dots,z_d, w_1,\dots, w_r)$ can be found 
such that
$$f(z,w)=w^2_1+\dots+w^2_r.$$

If $f$ has singularity of type A($d$) at $x$ then 
$S(f)=V(w_1,\dots,w_r)=J(f)$ so $j(f)=0$. There is a partial converse.

\lem 4 Suppose $f:\IC^{n+1},0\to \IC,0$. Suppose $I$ defines a complete 
intersection $\Sigma^d$ at $0$ with reduced structure, 
and suppose $j(f)=0$. Then $f$ has a singularity of type A($d$) at $0$.

\pf If $d=0$ the hypothesis implies that $J(f)=m_{n+1}$, and the result is 
implied by the  Morse lemma. Suppose $d>0$, then Theorem 5.14 p59 of 
\cite{P} implies that 
$\Sigma$ is an ICIS, and $f$ is A($d$) except perhaps at $0$. Further, the 
formula of 5.14 implies that the Tjurina number of $\Sigma$ is $0$, hence 
$\Sigma$ 
is smooth at the origin. Then proposition 3.13  p35 of \cite{P}, the 
formula cited above, and remark 5.3 on p52 imply that $f$ is A($d$) at the 
origin as well.

Now we are ready to prove our result about \Af.

\thm 5 Suppose $F:\IC^k\times\IC^{n+1},\IC^k\times 0,0\to \IC,0,0$, suppose the 
singular set of $F$, $S(F)$ is $V(I)$ where $I$ defines a 
family of complete intersections  with isolated singularities of fiber
dimension $d$, and every component of $V(I)$ contains $Y=\IC^k\times 0$. 
Suppose further that $J(F)=I$ off $Y$. Then:

1) If the pair $(\IC^k\times\IC^{n+1}-S(F),\IC^k\times 0)$ satisfies the 
\AF condition then $e(J(f_y),I_y,(y,0))$ is independent of $y$.

2) If $e(J(f_y),I_y,(y,0))$ is independent of $y$, then 
$\{\IC^k\times\IC^{n+1}-S(F), V(I)-Y, Y\}$ is an \AF stratification on 
some neighborhood of $Y$.

\pf To start the proof of 1), assume the \AF condition; this implies that  
${{\partial F} \over {\partial y_i}}\in \sdep
{J_z(F)}$ for $1\le i\le k$, by proposition 3.1. Now
$$e(J(F)(y),I(y),(y,z))=e(J_z(F)(y),I(y),(y,z))=e(J(f_y),I(y),(y,z))$$
\noindent for all (y,z) in some neighborhood of $(0,0)$. Since $J(F)=I$ 
off $Y$, this implies
$e(J(f_y),I(y),(y,z))=0$ off $Y$.

Since $\Gamma^k(I)=\Gamma^k(J_z(F))=\emptyset $, by the multiplicity-polar 
theorem, $$e(J(f_0),I(0),(0,0))=e(J(f_y),I(y),(y,0))$$ 
\noindent for all
$y$.

Now we prove 2). By hypothesis we have $I=J(F)$ off $Y$. So by lemma 3.4 
off of $Y$ we have that $V(I)$ is smooth and $F$ has only A($k+d$) 
singularities. So the pair
$\{\IC^k\times\IC^{n+1}-S(F), V(I)-Y\}$ has the \AF property.

Since $e(J(f_y),I_y,(y,0))$ is independent of $y$, and $I$ and $J_z(F)$ 
have no polars of dimension $k$, it then follows from the 
multiplicity-polar theorem
that $e(J(f_y),I_y,(y,z))=0$, for $z\ne 0$. This implies that 
$\?{J(f_y)}=\?{I_y}$. By lemma 3.3, $J(f_y)=I_y$. In turn this implies by 
lemma 3.5 that $V(I_y)$ 
is smooth off the origin and $f$ has an A($d$) singularity at points of 
$V(I_y)$ off the origin.  Now we have that $J_z(F)\subset I$ and at a 
point $(y,z)$ of $\Sigma$ off $Y$, 
$$\dim_{\IC} J(f_y,z)/(J(f_y,z)\cap m^2_z)=n+1-d\le \dim 
J_z(F)/(J_z(F)\cap m^2_{(y,z)})$$
$$\le \dim I/(I\cap m^2_{(y,z)})= n+1-d$$
Hence $J_z(F)=I$ at points of $\Sigma$ off $Y$.

Using what we have learned about $F$ above, we can describe the components 
of the exceptional divisor $E$ of 
$(B_{J_z(F)}(\IC^k\times\IC^{n+1}),\pi)$; 
we do this in order to get ready to apply 2) of 3.1, which will finish the 
proof. 

Let $\Sigma_i$ be the ith component of $\Sigma$; then there exists a 
component $V_i$ of $E$ which surjects to $\Sigma_i$. Suppose $V$ is a 
component of $E$
such that $\pi(V)$ is not contained in $Y$. Let $x$ be a point off $Y$ in 
$\pi(V)$. Then there is a neighborhood $U$ of $x$ in 
$\IC^k\times\IC^{n+1}$ such that on 
$U$,
$J(F)=J_z(F)=I$, and only one component of $\Sigma$ intersects $U$. Hence 
over $U$ the corresponding blowups are isomorphic; in particular there is 
only one
component of  each exceptional divisor which projects to $\Sigma\cap U$. 
So the $V_i$ are the only components of $E$ whose image does not lie in 
$Y$. Suppose $W$
is a component of $E$ whose image lies in $Y$. Then $W^{n+k}\subset 
Y^k\times \IP^{n}$, hence $W=Y^k\times \IP^{n}$ if $W$ exists. We have 
shown that every
component of $E$ projects to a set which  contains $Y$ in its closure. 
(This uses the fact that every $\Sigma_i$ contains $Y$ in its closure.)

Since \AF is true generically, there exists a Z-open set $U$ which 
contains a Z-open subset of $Y$, and on $U$ we have
 ${{\partial F} \over {\partial y_i}}\in \sdep
{J_z(F)}$ for $1\le i\le k$.  This implies that if we pull back $J_z(F)$ 
and $J(F)$ to the normalization of $B_{J_z(F)}(\IC^k\times\IC^{n+1})$, 
then along every component
of the exceptional divisor $E_N$ which meets $\pi_N^{-1}(U)$ in a Z-open 
set, that $\pi_N^*(J(F))=\pi_N^*(J_z(F)$. But this is true for all 
components of $E_N$,
since  every component of $E$ of
 $B_{J_z(F)}(\IC^k\times\IC^{n+1})$ projects to a set which 
contains $Y$ in its closure. This implies that $\?{J_z(F)}=\?{J(F)}$ at 
all points of $Y$ (\cite{LJT}). 

The last equality implies that $E_J$, the exceptional divisor of 
$B_{J(F)}(\IC^k\times\IC^{n+1})$, is finite over $E$. The components of 
$E_J$ which are in $\pi_N^{-1}(Y)$,
have dimension $k+n$ and have fiber dimension $n$, which is the fiber 
dimension of $W$, since they are finite over $W$. Hence they surject onto 
$W$, and hence $Y$. Since 
\AF holds generically, these components are in $C(Y)$, the conormal of 
$Y$, which also has dimension $n+k$, 
hence they are equal to the conormal, so there is only 1 such component.

Over each $V_i$ as we have seen there is only one component of $E_J$; 
since \AF holds between the open stratum and these components, 
a dimension count shows that this unique component is $C(\Sigma_i)$. The 
proof will be complete if we can show that each component of $\Sigma$
satisfies Whitney A over $Y$. (This is also what it means for \Af to hold 
for the pair $(\Sigma, Y)$.)

Claim: For every $i$, $C(\Sigma_i)\cap \pi_N^{-1}(Y\cap U)$ is dense in 
$C(\Sigma_i)\cap \pi_N^{-1}(Y)$. 

Since $C(\Sigma_i)\cap \pi_N^{-1}(Y\cap U)$ lies in $C(Y)$ this will 
finish the proof by 2) of 3.1.

By Lemma 5.7 p230 of \cite{GM}, we know that each component of 
$C(\Sigma_i)\cap \pi_N^{-1}(Y)$ has dimension $n+k-1$, that is, must be a 
hypersurface in $C(\Sigma_i)$. 
(This uses the fact that $I$ defines a complete intersection.) If the 
claim fails there must be a component for some $i$ of $C(\Sigma_i)\cap 
\pi_N^{-1}(Y)$ which does not surject onto
$Y$. Since $E_N|Y$ is finite over $E|Y\subset Y\times\IP^n$, this 
component must map to a subset of $Y$ of dimension $k-1$, and must have 
constant fiber dimension $n$.

Let $C$ be the fiber of the bad component over $0$. Consider 
$B_{J(F)(0)}(0\times\IC^{n+1})$. This must contain $C$ as a component of 
its exceptional divisor, as
$C$ is a subset of $B_{J(F)}(\IC^k\times\IC^{n+1})\cap 
0\times\IC^{n+1}\times\IP^{n+k}$, and its dimension is too small to be a 
component of the intersection. Construct a polar variety
of $J(F)$ of dimension $k+1$. This is a family of curves over $Y$; the 
fiber over $0$ contains a curve which is the projection of the 
intersection of 
the plane defining the polar with $B_{J(F)(0)}(0\times\IC^{n+1})$ forced 
by the existence of $C$. 
Let $\Gamma$ be the component of 
our polar which contains this curve.

We choose the plane of codimension $n$ of $\IP^{n+k}$ so that it misses 
the points of $C\cap C(Y)$. 
On some sufficiently small metric neighborhood of the origin in $\Gamma$, 
then we know that $\Gamma$ intersects $Y$ only at $(0,0)$. Restrict $I$ 
and $J(F)$ to $\Gamma$.
Now we apply the multiplicity-polar theorem again. $J(F)$ has no polar, 
because it is integrally dependent on $J_z(F)$ which has no polar. Over a 
generic $y$ value, 
the only points where $J(F)$ has support are on $\Sigma-Y$ hence 
$e(J(F)(y),I_y)=0$ at such points. We claim that the multiplicity of the 
pair 
$(J(F)(0),I_0)$ on $\Gamma(0)$ at $(0,0)$ is not zero. This number 
has an alternate meaning. It is part of the intersection number $ \int 
D_{J(F)(0),I(0)} \cdot l^{n}_{J(F)(0)}$, 
which in turn is part of $e((J(F)(0),I_0),(0,0))$ on $\IC^{n+1}$. We know 
that  $B_{\rho({J(F)(0)})}(\Projan {\cal R}(I_0))$, dominates both 
$B{_{I_0}}(\IC^{n+1})$ and 
$B{_{J(F)(0)}}(\IC^{n+1})$; corresponding to $C$ there is a component of 
the exceptional divisor of $B_{\rho({J(F)(0)})}(\Projan {\cal R}(I_0))$. 
The map to 
$B{_{I_0}}(\IC^{n+1})$
cannot be finite on this component, because the component projects to the 
origin in $\IC^{n+1}$, and the fiber dimension of the exceptional divisor 
of 
$B{_{I_0}}(\IC^{n+1})$ over  the origin must have dimension less than 
$n-d<n$, hence this component over $C$ makes a non-zero contribution to
$ \int D_{J(F)(0),I(0)} \cdot l^{n}_{J(F)(0)}$, so the multiplicity of the 
pair 
$(J(F)(0),I_0)$ on $\Gamma(0)$ at $(0,0)$ is not zero, so the 
multiplicity-polar theorem gives a contradiction--the change in 
multiplicity from the special fiber to the generic 
fiber is positive, but there is no polar variety of dimension $k$ of 
$J(F)$. So $C$ does not exist, which implies Whitney A holds for 
$(\Sigma-Y,Y)$  and the theorem is proved.

\rmk 6 The key point in the last proof, was the ability to take 
information about the $k+d$ dimensional strata of the total space, and 
relate it to the open stratum of $f_0$.
This was possible because we had good control on the conormals of the 
$k+d$ dimensional strata.

The above proof shows that it is easy to show that a stratification 
condition implies that the associated invariants are independent of 
parameter. To prove that the independence
from parameter implies the stratification condition requires in general 
the principle of specialization of integral dependence developed in 
\cite{G-7}.

As we shall see in general (remark 4.9) the  \AF condition does not imply 
that the L\^e numbers are independent of parameter.  We can introduce a
stronger notion of \AF which does imply that the L\^e numbers are constant 
if our ICIS is a curve. In the situation of Theorem 4.5 we say the strong 
\AF condition holds if
the \AF condition holds, and for a generic linear function $l$ the \Al 
condition holds for the pair $V(I)-Y,Y$

 From  Theorem 4.5, and the formula for $\lambda^0$ in Theorem 2.4, we can 
now show
that the
strong \AF condition implies that the L\^e numbers are constant in the 
setup originally  considered by Pellikaan.

\cor 7 Suppose $F:\IC^k\times\IC^{n+1},\IC^k\times 0,0\to \IC,0,0$, suppose the 
singular set of $F$, $S(F)$ is $V(I)$ where $I$ defines a 
family of complete intersection curves with isolated singularities, and 
every component of $V(I)$ contains $Y=\IC^k\times 0$. Suppose further that 
$J(F)=I$ off $Y$.
Suppose the pair $(\IC^k\times\IC^{n+1}-S(F),\IC^k\times 0)$ satisfies the 
\AF condition, and the pair $V(I)-\IC^k\times 0,\IC^k\times 0$ satisfies 
the \Al condition for a
generic linear function $l$, then the L\^e numbers of
$f_y$ at the origin are independent of $y$.

\pf Theorem 4.5 1) and Theorem 2.3 imply that $j(f_y)$ is constant along 
$Y$. 

The condition that the singular set of $F$ is $V(I)$ implies that $F$ is 
in $I^2$ (p8, prop 1.9 \cite{P}), hence $F=\sum\limits_{i,j}h_{i,j}g_ig_j$
where $\{g_i\}$ are a set of generators of $I$, and $h_{i,j}=h_{j,i}$ 
(\cite {P}, p54). Let $\Delta$ be the determinant of the matrix with 
entries $h_{i,j}$.
Then the number of D$_{\infty}$ points at $(y,z)$ is just the colength of 
$(\Delta_y)$ in $\O(V(I_y),z)$ (\cite{P}, p81 lemma 7.17). This number is 
just the local degree
at $(y,z)$ of the map with components $(\Delta, p)$ where $p$ is 
projection to the parameter space $Y$ on $V(I)$. 
Thus if $\delta(f_y)$ varies along $Y$ it must be upper semicontinuous, 
and if the  value for generic $y$ is less than the value over $y=0$,
 there must be other points in the fiber over $y$ where  $\delta(f_y)$ is 
non-zero.  However as the proof of Theorem 3.5 2) shows off $Y$ $f_y$ has 
only 
A$_{\infty}$ singularities on $V(I_y)$. Hence, $\delta(f_y)$ is constant 
along $Y$.

Since the pair $V(I)-\IC^k\times 0,\IC^k\times 0$ satisfies the \Al 
condition for a
generic linear function $l$, by Theorem 5.8 p232 of \cite{GM} the Milnor 
numbers of $V(I_y)$ and $V(I(y))\cap  V(l)$ are constant. Since $l$ is 
generic, the sum of these
Milnor numbers is just $e(JM(\Sigma_y)$, which is then independent of $y$. 
The result for $\lambda^0$ now follows from the formula for $\lambda^0$ in 
proposition 2.4.

Since the Milnor number of $V(I(y))\cap  V(l)$ is just the multiplicity of 
$\Sigma_y$, less $1$, the multiplicity of
$V(I_y)$ is independent of $Y$. Since the transverse Milnor number is 
always $1$, and 
the multiplicity of $V(I_y)$ constant, it follows that $\lambda^1$ is 
independent of $Y$ as well.

\cor 8 Suppose $F:\IC^k\times\IC^{n+1},\IC^k\times 0,0\to \IC,0,0$, suppose the 
singular set of $F$, $S(F)$ is $V(I)$ where $I$ defines a 
family of complete intersection curves with isolated singularities, and 
every component of $V(I)$ contains $Y=\IC^k\times 0$. Suppose further that 
$J(F)=I$ off $Y$.
Suppose the pairs $(\IC^k\times\IC^{n+1}-S(F),\IC^k\times 0)$,  
$V(I)-\IC^k\times 0,\IC^k\times 0$ satisfy the strong  \AF condition at 
$(0,0)$ then 

1) The homology of the Milnor fibre of of $f_y$ at the origin is 
independent of $y$ for all $y$ small.

If $n\ge 3$

2) The fibre homotopy-type of the Milnor fibrations of $f_y$ at the origin 
is independent of $y$ for all $y$ small.

If $n\ge 4$

3) The diffeomorphism-type of the Milnor fibrations of $f_y$ at the origin 
is independent of $y$ for all $y$ small.

\pf Since the strong \AF condition holds, Corollary 4.7 implies that the 
L\^e numbers are constant, then theorem 9.4 of \cite{Ma} p90 gives the 
result. (Although Massey 
states his theorem for the case where the dimension of the parameter 
stratum is 1, it also applies to the case at hand.)

This raises the interesting question of whether a strong \AF stratification or an \AF stratification
implies the triviality (in the sense of the last corollary) of the Milnor 
fibrations. The formulae in proposition 3.12 and remark 3.13 show that a strong \AF stratification implies that 
$\lambda^0(f_t)$ is independent of $t$ in these cases.

As the next example shows, in the \AF case,  this problem cannot be tackled by hoping that 
the existence of an \AF stratification implies that the L\^e numbers are 
constant. 

\rmk 9 This example shows that neither the \Af condition nor topological 
triviality imply that the L\^e numbers are constant. Let
$$f_t=z^5+ty^6z+y^7x+x^{15}.$$

This family of  functions was introduced by Brian\c con-Speder, 
(\cite{BS}) who  showed that $\mu^3(f_t):=\mu(f_t)=364$ for all $t$, while 
the Milnor number of a generic
hyperplane slice
$\mu^2(f_t)$ is $28$ when $t=0$ and $26$ otherwise. Historically, this example was 
important, because it showed that the $\mu_*$ constant condition was 
stronger than topological triviality.

Now consider $F_t=f_t^2+w^2$ where $w$ is a disjoint variable. Then

$$J_z(F)=<w,2f_t\Pd{f_t}x, 2f_t\Pd{f_t}y,2f_t\Pd{f_t}z>.$$

So the singular locus of $F$ is defined by $<w,f_t>$, hence is a family of 
complete intersections with isolated singularities.  A computation shows 
that:

$$e(J(F_t),(w,f_t))=j(F_t)=\mu^3(f_t).$$

Now, the only L\^e cycle of dimension 2 is $V(w,f_t)$, so 
$$\lambda^2(F_t)=m(X_t)=5,$$

while
$$\lambda^1(F_t)=m(\Gamma_1^1(X_t,0))=\mu^2(f_t)+\mu^1(f_t).$$

Now by 3.12

$$\lambda^0(F_t)=e(J(F_t),(w,f_t))+e(JM(V(w,f_t)))= \mu^3(f_t)+(\mu^2(f_t)+\mu^3(f_t)).$$

The first equality shows that the \AF condition holds by Theorem 3.5. 
However $\lambda^0(F_t)$ and $\lambda^1(F_t)$ vary with $t$. 
It is not hard to check by a vector field argument that the family of 
functions $F_t$ are topologically trivial; however this can be seen 
directly by the following argument
which was pointed out to me by J. N. Damon.

We know that there exists a topological trivialization $\phi(z, t) : 
\IC^{4}\to \IC^3$  of $f_t$, by \cite{BS}, so
		$f_t(\phi(z, t)) = f_0(z)$ 
Then, we can define 
	$	\Phi(z, w, t) =  (\phi(z, t), w) : \IC^{5} \to \IC^{4}$,
which gives a topological trivialization of $F_t$ since 
$$F_t(\Phi(x, w, t)) = f_t(\phi(x, t))  + w^2  = f_0(x) + w^2 =  F_0(x). 
$$

 Now we turn to the \Wf condition.  It is a paradox, but because this 
condition is stronger, it is easier to prove results about it.

\thm 10 Suppose $F:\IC^k\times\IC^{n+1},\IC^k\times 0,0\to \IC,0,0$, suppose the 
singular set of $F$, $S(F)$ is $V(I)$ where $I$ defines a 
family of complete intersections  with isolated singularities of fiber
dimension $d$, and every component of $V(I)$ contains $Y=\IC^k\times 0$. 
Suppose further that $J(F)=I$ off $Y$. Then:

1) If the pair $(\IC^k\times\IC^{n+1}-S(F),\IC^k\times 0)$ satisfies the 
\WF condition then $e(m_{n+1}J(f_y),I_y,(y,0))$ is independent of $y$.

2) If $e(m_{n+1}J(f_y),I_y,(y,0))$ is independent of $y$, then the pair 
$(\IC^k\times\IC^{n+1}-S(F),\IC^k\times 0)$ satisfies the \WF condition, 
and 
$\{\IC^k\times\IC^{n+1}-V(F), V(F)-V(I), V(I)-Y, Y\}$ is a Whitney 
stratification on some neighborhood of $Y$.

\pf 

\noindent 1) Suppose the pair $(\IC^k\times\IC^{n+1}-S(F),\IC^k\times 0)$ 
satisfies the \WF condition, then by Theorem 2.1 p23 of \cite{G-3}, 
the dimension of the fiber of the exceptional divisor over
$Y$ of $B_{m_YJ(F)}(\IC^k\times\IC^{n+1})$ is independent of $y$ and is 
$n$. This implies that the polar of dimension $k$ of $ m_YJ(F)$ is empty; 
hence by the multiplicity
polar theorem $e(m_{n+1}J(f_y),I_y,(y,0))$ is independent of $y$.

2) Suppose  $e(m_{n+1}J(f_y),I_y,(y,0))$ is independent of $y$. Off $Y$, 
$m_{n+1}J(f_y)=J(f_y)$, so off $Y$ by the same arguments found in the 
proof of 3.6, $J(f_y)=I_y$, so
 by the multiplicity polar theorem, $\Gamma^k(m_YJ(F))$ is empty, hence 
the dimension of the fiber of the exceptional divisor over
$Y$ of $B_{m_YJ(F)}(\IC^k\times\IC^{n+1})$ is $n$, hence is constant over 
$Y$. Then by Corollary 2.1, p 19 of \cite {G-3}, 
the pair $(\IC^k\times\IC^{n+1}-S(F),\IC^k\times 0)$ satisfies the \WF 
condition. This implies $V(F)-V(I)$ is Whitney over $Y$.

Since $F$ is of type A$_{\infty}$ off $Y$ it follows that $V(F)-V(I)$ is 
Whitney over $V(I)-Y$. It remains to show
$V(I)-Y$ is Whitney over $Y$. Suppose not; then for each $C$ and 
neighborhood $U$ of the origin there exists a sequence of points $x_i\in 
U$ on
 some component of $V(I)$, converging to the origin, and hyperplanes
$H_i$ which are tangent hyperplanes to $V(I)$ at $x_i$ such that

 $$\dist\bigl\leftp Y, H_{i}\bigr)> C\,\dist(x,Y).$$

From the proof of theorem 3.6, we have $C(V(I))\subset 
B_{J(F)}(\IC^k\times\IC^{n+1})$. This implies we can find points 
$\tilde x_i\in U\cap (\IC^k\times\IC^{n+1}-S(F))$ and hyperplanes $\tilde 
H_i$ tangent to the fibers of $F$ at $x_i$, such that
 the distance beteen $x_i$ and $\tilde x_i$, $H_i$ and $\tilde H_i$ is as 
small as desired. Then a similar inequality holds for $\tilde x_i$  and 
$\tilde H_i$, 
hence \WF fails, which is a contradiction.

\cor 11 Suppose in the above setup $e(m_{n+1}J(f_y),I_y,(y,0))$ is 
independent of $y$, then the family of functions $\{f_y\}$ is 
topologically trivial.

\pf Since $e(m_{n+1}J(f_y),I_y,(y,0))$ is independent of $y$, we have the 
pair $(\IC^k\times\IC^{n+1}-S(F),\IC^k\times 0)$ satisfies the \WF 
condition, and 
$\{\IC^k\times\IC^{n+1}-V(F), V(F)-V(I), V(I)-Y, Y\}$ is a Whitney 
stratification on some neighborhood of $Y$. Then we can lift the constant 
fields over $V(F)$, 
to the ambient space in such a way that the resulting fields can be 
integrated to give homeomorphisms.

There is a nice geometric interpretation of the number $e(m_{n+1}J(f_y))$ 
which we now describe. We denote the multiplicity of the relative polar 
variety of 
$f_y$ of dimension $i$ by $m^i(f_y)$.

\thm 12 Suppose $f:\IC^{n+1},0\to\IC,0$, $J$ any ideal in $\O_{n+1}$ such 
that $\dim_{\IC}J/J(f)<\infty$, then
$$e(m_{n+1}J(f),J)=e(J(f),J)+1+\sum\limits_{i=1}^{n}{{n+1}\choose i } 
m^i(f_y).$$

\pf This is exactly the content of the formula in Theorem 9.8 (i) p221 
\cite{K-T}.

\cor 13 Suppose $f:\IC^{n+1},0\to\IC,0$,  $S(f)$ is $V(I)$ where $I$ 
defines a 
 complete intersection  with isolated singularities of 
dimension $d$, and suppose further that $J(f)=I$ off $Y$. Then 
$$e(m_{n+1}J(f),I)=e(J(f),I)+1+\sum\limits_{i=1}^{n}{{n+1}\choose i } 
m^i(f_y).$$

\pf Follows immediately from Theorem 3.12

\cor 14 Suppose $F:\IC^k\times\IC^{n+1},\IC^k,0\to \IC,0,0$, suppose the 
singular set of $F$, $S(F)$ is $V(I)$ where $I$ defines a 
family of complete intersections  with isolated singularities of fiber
dimension $d$, and every component of $V(I)$ contains $Y=\IC^k\times 0$. 
Suppose further that $J(F)=I$ off $Y$. Then the following are equivalent:

1) $e(J(f_y),I_y)$ and the relative polar multiplicities of $f_y$ are 
independent of $y$.

2) \AF holds for the pair $(\IC^k\times\IC^{n+1}-V(I), Y)$, and the 
relative polar multiplicities of $f_y$ are independent of $y$.

3)  The pair $(\IC^k\times\IC^{n+1}-V(I),\IC^k\times 0)$ satisfies the \WF 
condition.

\pf 1) and 2) are equivalent by Theorem 3.5, while 2 and 3 are equivalent 
by Corollary 3.13 and  Theorem 3.9.

\sct References

\references

BMM
J. Brian\c con{,} P. Maisonobe and M. Merle,
 {\it Localisation de syst\`emes
diff\'erentiels, stratifications de Whitney et condition de Thom,}
 \invent 117 1994 531--50

BS
J. Brian\c con and J.P. Speder, {\it La trivialite topologique n'implique 
pas les conditions de Whitney}, C.R. Acad. Sc.Paris, v280 365-367 (1975)

B-R
 D. A. Buchsbaum and D. S. Rim,
  {\it A generalized Koszul complex. II. Depth and multiplicity,}
 \tams 111 1963 197--224

Bo1
J. Fern\'andez de Bobadilla, {\it Approximations of non-isolated 
singularities of finite codimension with respect to an isolated complete intersection 
singularity}.  Bull. London Math. Soc.  35  (2003),  no. 6, 812--816

F
 W. Fulton,
 ``Intersection Theory,''
 Ergebnisse der Mathematik und ihrer Grenzgebiete, 3. Folge
 $\cdot$ Band 2, Springer--Verlag, Berlin, 1984

 G-1
   T. Gaffney,
   {\it Aureoles and integral closure of modules,}
   in ``Stratifications, Singularities and Differential
 Equations~II'',  Travaux en Cours {\bf 55}, Herman, Paris, 1997, 55--62

 G-2
   T. Gaffney,
   {\it Integral closure of modules and Whitney equisingularity,}
   \invent 107 1992 301--22

 G-3
   T. Gaffney,
   {\it Plane sections, \Wf and \Af,} in ``Real and complex singularities 
(S\~ao Carlos, 1998)'',
 Chapman and Hall Res. Notes Math. {\bf 412}, 2000, 17-32

 G-4
   T. Gaffney,
   {\it Multiplicities and equisingularity of ICIS germs,}
    \invent 123 1996 209--220

G-5
	T. Gaffney, {\it Generalized Buchsbaum-Rim Multiplicities and a 
Theorem of Rees,} Communications in Algebra, vol 31 \#8 p3811-3828, 2003

G-6 
 T. Gaffney,{\it Polar methods, invariants of pairs of modules and 
equisingularity,}  Real and Complex Singularities (Sao Carlos, 2002),
Ed. T.Gaffney and M.Ruas, Contemp. Math.,\#354, Amer. Math. Soc.,
Providence,
RI, June 2004, 113-136

G-7
	T. Gaffney, {\it The Multiplicity-Polar Formula and Equisingularity}, 
in preparation

 G-G
   T. Gaffney and R. Gassler,
      {\it Segre numbers and hypersurface singularities,}
    J. Algebraic Geom. {\bf 8} (1999), 695--736

 GK
   T. Gaffney and S. Kleiman,
   {\it Specialization of integral dependence for modules}
   \invent 137 1999 541-574

 GK2
   T. Gaffney and S. Kleiman,
    {\it \Wf \hskip 2pt and specialization of integral dependence for 
modules,}
 in ``Real and complex singularities (S\~ao Carlos, 1998)'',
 Chapman and Hall Res. Notes Math. {\bf 412}, 2000, 33--45

 GM
   T. Gaffney and D. Massey,
   {\it Trends in equisingularity},  London Math. Soc. Lecture Note Ser.
 {\bf 263} (1999), 207--248

 Gas
   R. Gassler,
  {\it Segre numbers and hypersurface singularities}.
   Thesis, Northeastern University, 1999

H-M
 J.P.G. Henry and M. Merle,
 {\it Conormal Space and Jacobian module.  A short dictionary,}
 in ``Proceedings of the Lille Congress of Singularities," J.-P. Brasselet
(ed.), London Math. Soc. Lecture Notes {\bf 201} (1994), 147--74

J
		 G. Jiang,{\it  Functions with non-isolated singularities on 
singular spaces}, Thesis, Universiteit Utrecht, 1997

 K-T
   S. Kleiman and A. Thorup,
    {\it A geometric theory of the Buchsbaum--Rim multiplicity,}
   \ja 167 1994 168--231

 KT1
 S. Kleiman and A. Thorup,
  {\it The exceptional fiber of a generalized conormal space,}
  in ``New Developments in Singularity Theory." D.Siersma, C.T.C. Wall and 
V. Zakalyukin (eds.), Nato Science series,
II Mathematics, Physics and Chemistry-Vol. 21 2001 401-404

 Le-T
   D. T. L\^e and B. Teissier,
   {\it Cycles evanescents, sections planes et conditions de Whitney. II,}
 Proc. Sympos. Pure Math. {\bf 40}, part 2,
   Amer. Math.  Soc., 1983, 65--104

LJT
 M. Lejeune-Jalabert and B. Teissier,
 {\it Cl\^oture integrale des ideaux et equisingularit\'e, chapitre 1}  
Publ.
 Inst. Fourier   (1974)  

Ma
  D. Massey {\it L\^e Cycles and Hypersurface Singularities,}
  Springer Lecture Notes in Mathematics 1615 , (1995)

M
			D. Mumford, {\it Algebraic Geometry I Complex 
Projective Varieties,} Springer-Verlag 1976

N
A. Nemethi, {\it Hypersurface singularities with a 2 dimensional critical locus}, J. London Math. Soc. 
(3), 59, (1999), 922-938

P
		 R. Pellikaan,{\it  Hypersurface singularities and resolutions 
of Jacobi modules,} Thesis, Rijkuniversiteit Utrecht, 1985

P1
  R. Pellikaan, {\it Finite determinacy of functions with non-isolated 
singularities,} 
 Proc. London Math. Soc. vol. 57, pp. 1-26, 1988

S
		D. Siersma,{\it Isolated line singularities.}  Singularities, 
Part 2 (Arcata, Calif., 1981),  
		485--496, Proc. Sympos. Pure Math., 40, Amer. Math. Soc., 
Providence, RI, 1983

T-1
 B. Teissier,
 {\it The hunting of invariants in the geometry of the discriminant,}
 in ``Real and complex singularities, Oslo
1976,'' P. Holm (ed.), Sijthoff \& Noordhoff (1977), 565--678

 T-2
   B. Teissier,
   {\it Multiplicit\'es polaires, sections planes, et conditions de
 Whitney,}
   in ``Proc. La R\'abida, 1981.'' J. M. Aroca, R. Buchweitz, M. Giusti 
and
 M.  Merle (eds.), \splm 961 1982 314--491

Z 	
 A. Zaharia,{\it A study about singularities with non-isolated critical locus,} 
Thesis, Rijkuniversiteit Utrecht, 1993

Z1
		A. Zaharia, {\it Topological properties of certain 
singularities with critical 
locus a \break $2$-dimensional complete intersection,}  Topology Appl.  60  
(1994),  no. 2, 153--171

\endreferences

\bye